\def\author{ESTEVES, GAGN\'E, AND KLEIMAN}
\def\title{ABEL MAPS AND PRESENTATION SCHEMES}
\def\date{March 25, 2000}
 \def\abstract{%
 We sharpen the two main tools used to treat the
compactified Jacobian of a singular curve: Abel maps and presentation
schemes.  First we prove a smoothness theorem for bigraded Abel maps.
Second we study the two complementary filtrations provided by the images
of certain Abel maps and certain presentation schemes.  Third we study a
lifting of the Abel map of bidegree $(m,1)$ to the corresponding
presentation scheme.  Fourth we prove that, if a curve is blown up at a
double point, then the corresponding presentation scheme is a
$\IP^1$-bundle.  Finally, using Abel maps of bidegree $(m,1)$, we
characterize the curves having double points at worst.}
\def\PaperSize{letter} 

\def\TheMagstep{\magstep1}      
 \let\@=@  

\def\GetNext#1 {\def\NextOne{#1}\if\relax\NextOne\let\next=\relax
        \else\let\next=\DoIt \fi \next}
\def\DoIt{\Act\NextOne\GetNext}
\def\ActOn#1{\expandafter\GetNext #1\relax\ }
\def\defcs#1{\expandafter\xdef\csname#1\endcsname}

\parskip=0pt plus 1.75pt \parindent10pt
\hsize 32pc 
\vsize 48pc 
\abovedisplayskip 4pt plus3pt minus1pt
\belowdisplayskip=\abovedisplayskip
\abovedisplayshortskip 2.5pt plus2pt minus1pt
\belowdisplayshortskip=\abovedisplayskip

\def\TRUE{TRUE} 
\ifx\DoublepageOutput\TRUE \def\TheMagstep{\magstep0} \fi
\mag=\TheMagstep

\newskip\vadjustskip
\vadjustskip=0.5\normalbaselineskip
\def\centertext
 {\hoffset=\pgwidth \advance\hoffset-\hsize
  \advance\hoffset-2truein \divide\hoffset by 2\relax
  \voffset=\pgheight \advance\voffset-\vsize
  \advance\voffset-2truein \divide\voffset by 2\relax
  \advance\voffset\vadjustskip
 }
\newdimen\pgwidth\newdimen\pgheight
\def\letter{letter}\def\AFour{AFour}
\ifx\PaperSize\letter
 \pgwidth=8.5truein \pgheight=11truein
 \message{- Got a paper size of letter.  }\centertext
\fi
\ifx\PaperSize\AFour
 \pgwidth=210truemm \pgheight=297truemm
 \message{- Got a paper size of AFour.  }\centertext
\fi

\def\today{\ifcase\month\or     
 January\or February\or March\or April\or May\or June\or
 July\or August\or September\or October\or November\or December\fi
 \space\number\day, \number\year}
\nopagenumbers
 \newcount\pagenumber \pagenumber=1
 \def\advancepagenumber{\global\advance\pagenumber by 1}
\def\folio{\number\pagenum} 
\headline={%
  \ifnum\pagenum=0\hfill
  \else
   \ifnum\pagenum=1\firstheadline
   \else
     \ifodd\pagenum\oddheadline
     \else\evenheadline\fi
   \fi
  \fi
}
\expandafter\ifx\csname date\endcsname\relax \let\dato=\today
            \else\let\dato=\date\fi
\let\firstheadline\hfill
\def\oddheadline{\rm \eightpoint \rlap{\dato} \hfil
 \headtitle\hfil
 \llap{\folio}
 }
\def\evenheadline{\rm \eightpoint
 \rlap{\folio}
 \hfil\author
 \hfil\llap{\dato}
 }
\def\headtitle{\title}

 \newdimen\fullhsize \newbox\leftcolumn
 \def\fulline{\hbox to \fullhsize}
\def\doublepageoutput
{\let\lr=L
 \output={\if L\lr
           \global\setbox\leftcolumn=\columnbox \global\let\lr=R%
          \else \doubleformat \global\let\lr=L
          \fi
        \ifnum\outputpenalty>-20000 \else\dosupereject\fi
        }%
 \def\doubleformat{\shipout\vbox{%
     \ifx\PaperSize\AFour
           \fulline{\hfil\box\leftcolumn\hfil\columnbox\hfil}%
     \else
           \fulline{\hfil\hfil\box\leftcolumn\hfil\columnbox\hfil\hfil}%
     \fi             }%
     \advancepageno
}
 \def\columnbox{\vbox
   {\if E\topmark\headline={\hfil}\nopagenumbers\fi
    \makeheadline\pagebody\makefootline\advancepagenumber}%
   }%
\fullhsize=\pgheight \hoffset=-1truein
 \voffset=\pgwidth \advance\voffset-\vsize
  \advance\voffset-2truein \divide\voffset by 2
  \advance\voffset\vadjustskip
 
\ifx\FirstPageOnRight\TRUE 
 \null\vfill\nopagenumbers\eject\pagenum=1\relax
\fi
}
\ifx\DoublepageOutput\TRUE \let\pagenum=\pagenumber\doublepageoutput
 \else \let\pagenum=\pageno \fi

 \font\twelvebf=cmbx12          
 \font\smc=cmcsc10              
 \font\tenbi=cmmi14
 \font\sevenbi=cmmi10 \font\fivebi=cmmi7
 \newfam\bifam \def\bi{\fam\bifam} \textfont\bifam=\tenbi
 \scriptfont\bifam=\sevenbi \scriptscriptfont\bifam=\fivebi
 \mathchardef\variablemega="7121 \def\bigomega{{\bi\variablemega}}
 \mathchardef\variablenu="7117 
\catcode`\@=11          
\def\eightpoint{\eightpointfonts
 \setbox\strutbox\hbox{\vrule height7\p@ depth2\p@ width\z@}%
 \eightpointparameters\eightpointfamilies
 \normalbaselines\rm
 }
\def\eightpointparameters{%
 \normalbaselineskip9\p@
 \abovedisplayskip9\p@ plus2.4\p@ minus6.2\p@
 \belowdisplayskip9\p@ plus2.4\p@ minus6.2\p@
 \abovedisplayshortskip\z@ plus2.4\p@
 \belowdisplayshortskip5.6\p@ plus2.4\p@ minus3.2\p@
 }
\newfam\smcfam
\def\eightpointfonts{%
 \font\eightrm=cmr8 \font\sixrm=cmr6
 \font\eightbf=cmbx8 \font\sixbf=cmbx6
 \font\eightit=cmti8
 \font\eightsmc=cmcsc8
 \font\eighti=cmmi8 \font\sixi=cmmi6
 \font\eightsy=cmsy8 \font\sixsy=cmsy6
 \font\eightsl=cmsl8 \font\eighttt=cmtt8}
\def\eightpointfamilies{%
 \textfont\z@\eightrm \scriptfont\z@\sixrm  \scriptscriptfont\z@\fiverm
 \textfont\@ne\eighti \scriptfont\@ne\sixi  \scriptscriptfont\@ne\fivei
 \textfont\tw@\eightsy \scriptfont\tw@\sixsy \scriptscriptfont\tw@\fivesy
 \textfont\thr@@\tenex \scriptfont\thr@@\tenex\scriptscriptfont\thr@@\tenex
 \textfont\itfam\eightit        \def\it{\fam\itfam\eightit}%
 \textfont\slfam\eightsl        \def\sl{\fam\slfam\eightsl}%
 \textfont\ttfam\eighttt        \def\tt{\fam\ttfam\eighttt}%
 \textfont\smcfam\eightsmc      \def\smc{\fam\smcfam\eightsmc}%
 \textfont\bffam\eightbf \scriptfont\bffam\sixbf
   \scriptscriptfont\bffam\fivebf       \def\bf{\fam\bffam\eightbf}%
 \def\rm{\fam0\eightrm}%
 }
\def\vfootnote#1{\insert\footins\bgroup
 \eightpoint\catcode`\^^M=5\leftskip=0pt\rightskip=\leftskip
 \interlinepenalty\interfootnotelinepenalty
  \splittopskip\ht\strutbox 
  \splitmaxdepth\dp\strutbox \floatingpenalty\@MM
  \leftskip\z@skip \rightskip\z@skip \spaceskip\z@skip \xspaceskip\z@skip
  \textindent{#1}\footstrut\futurelet\next\fo@t}

\def\p.{p.\penalty\@M \thinspace}
\def\pp.{pp.\penalty\@M \thinspace}
\newcount\sctno \newskip\sctnskip \sctnskip=0pt plus\baselineskip
\def\sctn#1\par
  {\removelastskip\vskip\sctnskip 
  \vskip-\sctnskip \bigskip\medskip
  \centerline{#1}\nobreak\medskip
}

\def\sct#1 {\sctno=#1\relax\sctn#1. }

\def\item#1 {\par\indent\indent\indent
 \hangindent3\parindent
 \llap{\rm (#1)\enspace}\ignorespaces}
 \def\inpart#1 {{\rm (#1)\enspace}\ignorespaces}
 \def\part {\par\inpart}

\def\Cs#1){\(\number\sctno.#1)}
\def\part#1 {\par\(#1)\enspace\ignorespaces}

\def\dsc#1 #2.{\medbreak\Cs#1) {\it #2.} \ignorespaces}
\def\proclaim#1 #2 {\medbreak
  {\bf#1 (\number\sctno.#2).  }\ignorespaces\bgroup\it}
\def\endproclaim{\par\egroup\medskip}
\def\pf{\endproclaim{\bf Proof.} \ignorespaces}
\def\lem{\proclaim Lemma } \def\prp{\proclaim Proposition }
\def\cor{\proclaim Corollary }  \def\thm{\proclaim Theorem }
\def\dfn#1 {\medbreak {\bf Definition (\number\sctno.#1).  }\ignorespaces}
\def\rmk#1 {\medbreak {\bf Remark (\number\sctno.#1).  }\ignorespaces}
\def\eg#1 {\medbreak {\bf Example (\number\sctno.#1).  }\ignorespaces}
 \newcount\refno \refno=0        \def\NoKey{*!*}
 \def\MakeKey{\advance\refno by 1 \expandafter\xdef
  \csname\TheKey\endcsname{{\number\refno}}\NextKey}
 \def\NextKey#1 {\def\TheKey{#1}\ifx\TheKey\NoKey\let\next\relax
  \else\let\next\MakeKey \fi \next}
 \def\RefKeys #1\endRefKeys{\expandafter\NextKey #1 *!* }
 \def\SetRef#1 #2,{\hang\llap
  {[\csname#1\endcsname]\enspace}{\smc #2},}
 \newbox\keybox \setbox\keybox=\hbox{[25]\enspace}
 \newdimen\keyindent \keyindent=\wd\keybox
\def\references{\kern-\medskipamount
  \sctn References\par
  \vskip-\medskipamount
  \bgroup   \frenchspacing   \eightpoint
   \parindent=\keyindent  \parskip=\smallskipamount
   \everypar={\SetRef}\par}
\def\endreferences{\egroup}

 \def\serial#1#2{\expandafter\def\csname#1\endcsname ##1 ##2 ##3
        {\unskip\ {\it #2\/} {\bf##1} (##2), ##3.}} 

\def\uts{\penalty\@M \thinspace\ignorespaces}
\def\(#1){{\let~=\uts\rm(#1)}}
\def\relaxnext@{\let\next\relax}
\def\cite#1{\relaxnext@
 \def\nextiii@##1,##2\end@{\unskip\space{\rm[\SetKey{##1},\let~=\uts##2]}}%
 \in@,{#1}\ifin@\def\next{\nextiii@#1\end@}\else
 \def\next{{\rm[\SetKey{#1}]}}\fi\next}
\newif\ifin@
\def\in@#1#2{\def\in@@##1#1##2##3\in@@
 {\ifx\in@##2\in@false\else\in@true\fi}%
 \in@@#2#1\in@\in@@}
\def\SetKey#1{{\bf\csname#1\endcsname}}

\catcode`\@=12  

\let\:=\colon \let\ox=\otimes \let\x=\times
\let\wt=\widetilde \let\wh=\widehat
\def\dr#1{#1^\dagger}
\let\into=\hookrightarrow \def\onto{\to\mathrel{\mkern-15mu}\to}
 
\def\smashedlongrightarrow{\setbox0=\hbox{$\longrightarrow$}\ht0=1pt\box0}
\def\risom{\buildrel\sim\over{\smashedlongrightarrow}}
 \def\lgto{-\mathrel{\mkern-10mu}\to}
 \def\smashedlgto{\setbox0=\hbox{$\scriptstyle\lgto$}\ht0=1.85pt
        \lower1.25pt\box0}
\def\tto{\buildrel\lgto\over{\smashedlgto}}
\let\vf=\varphi  \let\?=\overline

\def\CJ{{\bar J}}  

 \def\IP{{\bf P}}  \def\Im{{\bf m}}
\def\Act#1{\defcs{c#1}{{\cal#1}}}
 \ActOn{A B C D E F G H I J K L M N O P Q R T }
\def\Act#1{\defcs{#1}{\mathop{\rm#1}\nolimits}}
 \ActOn{cod Div Ext Hilb length colength Pic Quot sm Spec Supp Grass Hom }
\def\Act#1{\defcs{c#1}{\mathop{\it#1}\nolimits}}
 \ActOn{Cok Ext Hom Ker Sym Im }
\def\Act#1{\defcs{#1}{\hbox{\rm\ #1 }}}
 \ActOn{and by for where with on }

\catcode`\@=11

 \def\activeat#1{\csname @#1\endcsname}
 \def\def@#1{\expandafter\def\csname @#1\endcsname}
 {\catcode`\@=\active \gdef@{\activeat}}

\let\ssize\scriptstyle
\newdimen\ex@   \ex@.2326ex

 \def\requalfill{\cleaders\hbox{$\mkern-2mu\mathord=\mkern-2mu$}\hfill
  \mkern-6mu\mathord=$}
 \def\eqfill{$\m@th\mathord=\mkern-6mu\requalfill}
 \def\deffill{\hbox{$:=$}$\m@th\mkern-6mu\requalfill}
 \def\fiberbox{\hbox{$\vcenter{\hrule\hbox{\vrule\kern1ex
     \vbox{\kern1.2ex}\vrule}\hrule}$}}

 \font\arrfont=line10
 \def\Swarrow{\vcenter{\hbox{$\swarrow$\kern-.26ex
    \raise1.5ex\hbox{\arrfont\char'000}}}}

 \newdimen\arrwd
 \newdimen\minCDarrwd \minCDarrwd=2.5pc
        
 \def\findarrwd#1#2#3{\arrwd=#3%
  \setbox\z@\hbox{$\ssize\;{#1}\;\;$}%
 \setbox\@ne\hbox{$\ssize\;{#2}\;\;$}%
  \ifdim\wd\z@>\arrwd \arrwd=\wd\z@\fi
  \ifdim\wd\@ne>\arrwd \arrwd=\wd\@ne\fi}
 \newdimen\arrowsp\arrowsp=0.375em
 \def\findCDarrwd#1#2{\findarrwd{#1}{#2}{\minCDarrwd}
    \advance\arrwd by 2\arrowsp}
 \newdimen\minarrwd 
 \setbox\z@\hbox{$\longrightarrow$} \minarrwd=\wd\z@

 \def\harrow#1#2#3#4{{\minarrwd=#1\minarrwd%
   \findarrwd{#2}{#3}{\minarrwd}\kern\arrowsp
    \mathrel{\mathop{\hbox to\arrwd{#4}}\limits^{#2}_{#3}}\kern\arrowsp}}
 \def@]#1>#2>#3>{\harrow{#1}{#2}{#3}\rightarrowfill}
 \def@>#1>#2>{\harrow1{#1}{#2}\rightarrowfill}
 \def@<#1<#2<{\harrow1{#1}{#2}\leftarrowfill}
 \def@={\harrow1{}{}\eqfill}
 \def@:#1={\harrow1{}{}\deffill}
 \def@ N#1N#2N{\vCDarrow{#1}{#2}\UpDownarrow}
 \def\UpDownarrow{\uparrow\,\Big\downarrow}

\def\hookleftarrowfill{$\m@th\mathord\leftarrow\mkern-6mu%
  \cleaders\hbox{$\mkern-2mu\mathord-\mkern-2mu$}\hfill\mkern-12mu%
  \relbar\joinrel\rhook$}
 \def@ l#1<#2<{\harrow1{#1}{#2}\hookleftarrowfill}

\def@'#1'#2'{\harrow1{#1}{#2}\tarrowfill}
 \def\lgTo{\dimen0=\arrwd \advance\dimen0-2\arrowsp
        \hbox to\dimen0{\rightarrowfill}}
 \def\smashedlgTo{\setbox0=\hbox{$\scriptstyle\lgTo$}\ht0=1.85pt
        \lower1.25pt\box0}
 \def\tto{\buildrel\lgTo\over{\smashedlgTo}}
 \def\tarrowfill{\hfil$\tto$\hfil}  

 \def@={\ifodd\row\harrow1{}{}\eqfill
   \else\vCDarrow{}{}\Vert\fi}
 \def@.{\ifodd\row\relax\harrow1{}{}\hfill
   \else\vCDarrow{}{}.\fi}
 \def@|{\vCDarrow{}{}\Vert}
 \def@ V#1V#2V{\vCDarrow{#1}{#2}\downarrow}
\def@ A#1A#2A{\vCDarrow{#1}{#2}\uparrow}
 \def@(#1){\arrwd=\csname col\the\col\endcsname\relax
   \hbox to 0pt{\hbox to \arrwd{\hss$\vcenter{\hbox{$#1$}}$\hss}\hss}}

 \def\squash#1{\setbox\z@=\hbox{$#1$}\finsm@@sh}
\def\finsm@@sh{\ifnum\row>1\ht\z@\z@\fi \dp\z@\z@ \box\z@}

 \newcount\row \newcount\col \newcount\numcol \newcount\arrspan
 \newdimen\vrtxhalfwd  \newbox\tempbox

 \def\innernewdimen{\alloc@1\dimen\dimendef\insc@unt}
 \def\measureinit{\col=1\vrtxhalfwd=0pt\arrspan=1\arrwd=0pt
   \setbox\tempbox=\hbox\bgroup$}
 \def\setinit{\col=1\hbox\bgroup$\ifodd\row
   \kern\csname col1\endcsname
   \kern-\csname row\the\row col1\endcsname\fi}
 \def\findvrtxhalfsum{$\egroup
  \expandafter\innernewdimen\csname row\the\row col\the\col\endcsname
  \global\csname row\the\row col\the\col\endcsname=\vrtxhalfwd
  \vrtxhalfwd=0.5\wd\tempbox
  \global\advance\csname row\the\row col\the\col\endcsname by \vrtxhalfwd
  \advance\arrwd by \csname row\the\row col\the\col\endcsname
  \divide\arrwd by \arrspan
  \loop\ifnum\col>\numcol \numcol=\col%
 \expandafter\innernewdimen \csname col\the\col\endcsname
     \global\csname col\the\col\endcsname=\arrwd
   \else \ifdim\arrwd >\csname col\the\col\endcsname
      \global\csname col\the\col\endcsname=\arrwd\fi\fi
   \advance\arrspan by -1 %
   \ifnum\arrspan>0 \repeat}
 \def\setCDarrow#1#2#3#4{\advance\col by 1 \arrspan=#1
    \arrwd= -\csname row\the\row col\the\col\endcsname\relax
    \loop\advance\arrwd by \csname col\the\col\endcsname
     \ifnum\arrspan>1 \advance\col by 1 \advance\arrspan by -1%
     \repeat
    \squash{\mathop{
     \hbox to\arrwd{\kern\arrowsp#4\kern\arrowsp}}\limits^{#2}_{#3}}}
 \def\measureCDarrow#1#2#3#4{\findvrtxhalfsum\advance\col by 1%
   \arrspan=#1\findCDarrwd{#2}{#3}%
    \setbox\tempbox=\hbox\bgroup$}
 \def\vCDarrow#1#2#3{\kern\csname col\the\col\endcsname
    \hbox to 0pt{\hss$\vcenter{\llap{$\ssize#1$}}%
     \Big#3\vcenter{\rlap{$\ssize#2$}}$\hss}\advance\col by 1}

 \def\setCD{\def\harrow{\setCDarrow}%
  \def\\{$\egroup\advance\row by 1\setinit}
  \m@th\lineskip3\ex@\lineskiplimit3\ex@ \row=1\setinit}
 \def\endsetCD{$\strut\egroup}
 \def\measure{\bgroup
  \def\harrow{\measureCDarrow}%
  \def\\##1\\{\findvrtxhalfsum\advance\row by 2 \measureinit}%
  \row=1\numcol=0\measureinit}
 \def\endmeasure{\findvrtxhalfsum\egroup}

\newbox\CDbox \newdimen\sdim

 \newcount\savedcount
 \def\CD#1\endCD{\savedcount=\count11%
   \measure#1\endmeasure
   \vcenter{\setCD#1\endsetCD}%
   \global\count11=\savedcount}
 \catcode`\@=\active


 \RefKeys
 AIK76 AK76 AK79II AK80 AK90 E95 EGK99 EGA KK81 R80 ZS60
 \endRefKeys
\font\ssi=cmssi10
\def\firstheadline{\hfil\ssi
 To Robin, who taught us geometry from Euclid to Grothendieck}
{\leftskip=0pt plus1fill \rightskip=\leftskip
 \obeylines
 \null\bigskip\bigskip\bigskip\bigskip
 {\twelvebf \title
 } \bigskip
 \footnote{}{\noindent %
 MSC-class: 14H40 (Primary) 14H20 (Secondary).}
 Eduardo Esteves\footnote{$^{1}$}{%
    Supported in part by PRONEX and CNPq Proc. 300004/95-8 (NV).}
 {\eightpoint\it\medskip
 Instituto de Matem\'atica Pura e Aplicada
 Estrada D. Castorina {\sl110, 22460--320} Rio de Janeiro RJ, BRAZIL
 \rm E-mail: \tt Esteves\@impa.br \medskip
 }
 Mathieu Gagn\'e
 {\eightpoint\it\medskip
 EMC Corporation
 171 South Street, Hopkinton, MA {\sl01748}, USA
  \rm E-mail: \tt MGagne\@emc.com \medskip
 } and \medskip
 Steven Kleiman\footnote{$^{2}$}{%
    Supported in part by the US NSF, the Brazilian CNPq, and the Danish
NSRC.}
 {\eightpoint\it\medskip
 Department of Mathematics, Room {\sl 2-278} MIT,
 {\sl77} Mass Ave, Cambridge, MA {\sl02139-4307}, USA
 \rm E-mail: \tt Kleiman\@math.mit.edu \medskip
 \rm \dato \bigskip
 }
}
{\parindent=1.5\parindent \narrower \eightpoint \noindent
 {\smc Abstract.}\enspace \ignorespaces \abstract \par}

\sct1 Introduction

Let $C$ be an integral projective curve over a field.  For the most
part, we'll work with flat families of such $C$ over an arbitrary
locally Noetherian base scheme.  However, for the sake of simplicity in
this introduction, we'll fix $C$ and assume that the base field is
algebraically closed of arbitrary characteristic.  Given $n$, let
$J^n_C$ denote the component of the Picard scheme parameterizing
invertible sheaves of degree $n$, and let $\CJ^n_C$ denote its natural
compactification by torsion-free rank-1 sheaves.  These schemes are
(essentially) independent of $n$, and they are often called the
(generalized) {\it Jacobian\/} of $C$ and its {\it compactified
Jacobian}.

To treat $\CJ^n_C$, there are two main tools: Abel maps and presentation
schemes.  Abel maps were studied in \cite{AK80}, and presentation
schemes in \cite{AK90}.  In the present paper, we'll advance those
studies, and tie them together.
  Notably, we'll develop the theory we
need to be able in \cite{EGK99} to prove the following {\it autoduality
theorem:\/} given an invertible sheaf $\cI$ of degree 1 on $C$, the
corresponding Abel map $A_{\cI,C}\:C\to\CJ^0_C$ induces an isomorphism,
        $$A_{\cI,C}^*\:\Pic_{\CJ_C^0}^0\risom J_C^0,$$
 and $A_{\cI,C}^*$ is independent of $\cI$, at least when $C$ has at
worst double points.  (We'll discuss the Abel map next, and we'll
outline the proof of the autoduality theorem at the end of the
introduction.)  However, in the present paper, we'll develop the theory
in its natural generality, not simply in the special cases we'll need in
\cite{EGK99}; moreover, this general theory is of independent interest.

Abel maps are forgetful maps.  Let $\cF$ be a torsion-free rank-1 sheaf
of degree $f$ on $C$, and $\cI$ an invertible sheaf of degree $m$.  By
definition, the {\it Abel\/} map,
  $$A_{\cI,\cF}\:\Quot^n_\cF\to\CJ^{m+f-n}_C,$$
 carries a subsheaf $\cJ\subset\cF$ to the isomorphism class of
$\cI\ox\cJ$, forgetting $\cF$.  In practice, two particular choices of
$\cF$ are especially important: (1) $\cF=\cO_C$, whence $\Quot^1_\cF=C$
and $\Quot^n_\cF= \Hilb^n_C$, and (2) $\cF=\bigomega$, where $\bigomega$
is the dualizing sheaf.  As above, when $\cF=\cO_C$, we'll write simply
$A_{\cI,C}$.

We'll study the $A_{\cI,\cF}$ by bundling them together, as $\cI$
varies but $\cF$ is fixed, to form a {\it bigraded Abel map} of {\it
bidegree} $(m,n)$,
  $$A_\cF\:J_C^m\x \Quot^n_\cF\to\CJ^{m+f-n}_C.$$
  When $\cF=\cO_C$, we'll write $A_{C}$.  Now, $A_\cF$ extends
naturally over the subset of $\CJ_C^m\x \Quot^n_\cF$ of pairs
$(\cI,\cJ)$ such that $\cI$ is invertible where $\cJ$ is not.  By
Corollary (2.6), the extension of $A_C$ is smooth and surjective if
$C$ is Gorenstein.  In \cite{EGK99}, we'll only use $A_C$ on $J_C^1\x
C$, but we'll make essential use of its smoothness, or rather, of its
flatness.

Corollary (2.6) follows from Theorem (2.4), which can be formulated as
follows.  For any sheaf $\cG$ on $C$, set
$\cG^\star:=\cHom_C(\cG,\bigomega)$.  For any $m$ and $n\ge0$, let
$W^{m,n}_\cF$ be the subset of $\CJ_C^m\x \Quot^n_\cF$ of pairs
$(\cI,\cJ)$ such that $\cI^\star$ is invertible where $\cJ^\star$ is
not.  Let $g$ be the arithmetic genus of $C$, and define another {\it
bigraded Abel map} of {\it bidegree} $(m,n)$,
         $$\alpha_\cF\:W^{m,n}_\cF \to \CJ_C^{f-2g+2+m-n}
    \hbox{ by }\alpha_\cF(\cI,\cJ):=(\cI^\star\ox\cJ^\star)^\star.$$
 Theorem~(2.4) asserts that $\alpha_\bigomega$ is smooth and
 surjective.  We'll prove it by checking the infinitesimal criterion,
which is usually the most manageable way to prove smoothness, especially
when the source and target are not smooth.

  Since $\alpha_\bigomega$ is smooth, it's flat, so open.  Hence, we can
construct an open filtration of $\CJ^n_C$ as follows.  Fix any closed
subset $Y$ of $C$.  Let $\CJ^m_{C,Y}\subset \CJ^m_C$ be the open subset
of sheaves $\cI$ such that $\cI^\star$ is invertible along $Y$.  Set
     $$W^{m,k}_{Y,\cF}:= W^{m,k}_\cF\cap (\CJ_{C,Y}^m\x\Quot^k_\cF)
 \and\CJ^{n,k}_{C,Y}:=\alpha_{\bigomega}(W^{n+k,k}_{Y,\bigomega}).$$
 Let $\delta$ denote the drop in the genus from $C$ to the partial
normalization along $Y$.  Then, by Proposition~(4.4), the open subsets
$\CJ^{n,k}_{C,Y}$ provide an exhaustive increasing filtration,
        $$\CJ^n_{C,Y}=\CJ^{n,0}_{C,Y}\subseteq\dots
        \subseteq\CJ^{n,\delta}_{C,Y}=\CJ^{n,\delta+1}_{C,Y}
        =\cdots=\CJ^n_C.$$

For example, take $Y=C$.  Then $\delta$ is the difference between the
arithmetic genus and the geometric genus of $C$.  Now, for every $m$,
define a map
    $$\iota\: \CJ^m_C\to\CJ^{2g-2-m}_C\hbox{ by }\iota(\cI):=\cI^\star.$$
 Then $\iota^2=1$; so $\iota$ is an isomorphism, and $\iota^{-1}=\iota$.
We have $\CJ^m_{C,C}=\iota(J^{2g-2-m}_C)$.  Hence, $\CJ^m_{C,C}$ is
smooth.  Moreover, $W^{m,k}_{C,\bigomega}$ is equal to
$\CJ^m_{C,C}\x\Quot^k_\bigomega$.  Since $\alpha_\bigomega$ is smooth,
it follows that, locally in the smooth topology, $\CJ^{n,k}_{C,C}$ is
isomorphic to $\Quot^k_\bigomega$; hence, these two schemes have
singularities of the same type.  Since $\CJ^{n,\delta}_{C,C}= \CJ^n_C$,
the singularities of $\CJ^n_C$ are thus of the same type as those of
$\Quot^k_\bigomega$ for $k\ge\delta$.  Previously, this result was known
 only for $k\ge 2g-1$; in this case, it is a consequence of the
 smoothness and surjectivity of $A_{\cO_C,\bigomega}$, which was proved
in \cite{AK80, (8.4), \p.103}.

If $C$ is Gorenstein, then tensoring with $\bigomega$ defines an
isomorphism $\Hilb^n_C\risom\Quot^n_\bigomega$; in particular, then
$C\risom\Quot^1_\bigomega$.  Therefore, the singularities of
$\CJ^{n,1}_{C,C}$ are then of the same type as those of $C$.
 Furthermore, Theorem (6.7) asserts that, $C$ being Gorenstein, $C$ has
at worst double points if and only if $\CJ^{n,1}_{C,C}$ is large enough
to contain all the points of $\CJ^n_C$ of dimension at least $g-1$.  The
latter condition may be replaced by the bound $\dim(\CJ^n_C-V)\le g-2$
where $V$ stands for the image of
$J_C^{n+1}\x C$
 under $A_C$, according to Corollary~(6.8).  Indeed, since $C$ is
Gorenstein, a torsion-free rank-1 sheaf $\cI$ on $C$ is invertible if
and only if $\cI^\star$ is invertible; whence
$\CJ^{n+1}_{C,C}=J^{n+1}_C$, and so $\CJ^{n,1}_{C,C}=V$.

 Presentation schemes arise when a birational map $\vf\:\dr C\to C$ is
given.  Although $\vf$ induces a map on the Jacobians $\vf^*\: J^n_C\to
J^n_{\dr C}$, nevertheless, when $\vf$ is nontrivial, $\vf^*$ does not
extend to a regular map from $\CJ^n_C$ to $\CJ^n_{\dr C}$ (because such
an extension would restrict to an inverse of $\vf$, where $C$ and $\dr
C$ are viewed as embedded by $A_{\cI,C}$ and $A_{\vf^*\cI,\dr C}$ for
any given invertible sheaf $\cI$ on $C$ of degree $n+1$, see \cite{AK80,
(8.8), \p.108}).  To obtain a modification of $\CJ^n_C$ on which $\vf^*$
does extend, we could form the graph of $\vf^*$, and take its closure in
$\CJ^n_C\x\CJ^n_{\dr C}$.  However, there is a more useful modification,
and it is one of a series of projective moduli spaces $P^{n,r}_\vf$,
known as the ``presentation schemes.''

The {\it presentation scheme} $P^{n,r}_\vf$ is defined as follows; see
(3.1).  Let $Y\subset C$ be the subset where $\vf\:\dr C\to C$ fails to
be an isomorphism, and give $Y$ the subscheme structure defined by the
conductor ideal $\cI_Y$ of $\vf$.
Then $P^{n,r}_\vf$ parameterizes the pairs $(\cI,\dr\cI)$ consisting
of a torsion-free rank-1 sheaf $\dr\cI$ of degree $n$ on $C^\dagger$ and
a subsheaf $\cI$ of $\vf_*\dr\cI$ such that $\vf_*\dr\cI/\cI$ is an
$\cO_Y$-module of length $r$; in other words, $P^{n,r}_\vf$
parameterizes {\it presentations}, which are the corresponding
injective maps $h\:\cI \to \vf_*\dr\cI$.

The presentation scheme $P^{n,r}_\vf$ is useful because there are
natural maps,
        $$\kappa_\vf\:P^{n,r}_\vf \to \CJ^{n+\delta-r}_C \and
        \pi_\vf\:P^{n,r}_\vf \to \CJ^n_{C^\dagger}, $$
 where $\delta$ is the genus drop from $C$ to $\dr C$.
These maps
 are defined simply by sending a pair $(\cI,\dr\cI)$ to $\cI$ and
to $\dr\cI$.
 If we understand $\pi_\vf$ and $\kappa_\vf$ well enough, then we may
now derive properties of $\CJ^{n+\delta-r}_C$ from those of
$\CJ^n_{C^\dagger}$.

In particular, $P^{n,\delta}_\vf$ is
 a suitable
 modification of
$\CJ^n_C$.  Indeed, restricting $\kappa_\vf$ yields an isomorphism,
        $$\kappa_\vf^{-1}(J^n_C)\risom J^n_C;$$
 its inverse is defined by sending an invertible sheaf $\cL$ on $C$ to
the canonical presentation $\cL\to\vf_*\vf^*\cL$.  Hence, $\pi_\vf$
extends $\vf^*$.

When $C$ and $\dr C$ are Gorenstein, Corollary (5.5) describes
$\kappa_\vf$, locally in the smooth topology, over a larger open subset
of $\CJ^n_C$ than $J^n_C$.  More precisely, Corollary (5.5) says that
there is a map $\Lambda_\vf$ that completes the product square in the
following diagram:
 $$ \CD
 \CJ^{n+1}_C\x\dr C @l<< (1\x\vf)^{-1}W^{n+1,1}_{Y,C}
        @>\Lambda_\vf>> P^{n,\delta}_\vf@>\pi_\vf>> \CJ^n_{\dr C}\\
      @V1\x\vf VV  @VVV  @(\fiberbox) @V\kappa_\vf VV \\
 \CJ^{n+1}_C\x C @l<< W^{n+1,1}_{Y,C}
        @>A_C>> \CJ^n_C
 \endCD$$
 where $W^{n+1,1}_{Y,C}:=W^{n+1,1}_{Y,\cO_C}$.
 Note that, since $C$ is Gorenstein, $A_C$ extends over
$W^{n+1,1}_{Y,C}$, and $A_C(W^{n+1,1}_{Y,C})$ contains $J^n_C$.

 Corollary (5.5) follows from Theorem (5.4), which makes a similar
statement about $\alpha_\cF\:W^{m,1}_{Y,\cF}\to \CJ^r_C$ where
$r:=m-1+f-2g+2$ for any $\cF$ of degree $f$.  On the other hand,
Corollary (5.5) generalizes \cite{AK90, (17)}, which makes a cruder
assertion in the
 special
 case where $Y$ is an ordinary node or cusp.

When $C$ and $\dr C$ are Gorenstein and also $\dr C$ is smooth along
$\dr Y:=\vf^{-1}Y$, Theorem (4.9) describes, in terms of $\kappa_\vf$,
the complementary decreasing closed filtration,
        $$(\CJ^n_C-\CJ^{n,0}_{C,Y})\supseteq(\CJ^n_C-\CJ^{n,1}_{C,Y})
  \supseteq\cdots\supseteq(\CJ^n_C-\CJ^{n,\delta}_{C,Y})=\emptyset.$$
 Namely, Theorem (4.9) asserts that
   $$\CJ^n_C-\CJ^{n,k}_{C,Y}=\iota\kappa_\vf(P^{m,\delta-k-1}_\vf)
        \hbox{ where }m:=2g-3-k-n;$$
 moreover, $\iota\kappa_\vf$ restricts to an isomorphism of schemes,
        $$(\iota\kappa_\vf)^{-1}(\CJ^{n,k+1}_{C,Y}-\CJ^{n,k}_{C,Y})
        \risom(\CJ^{n,k+1}_{C,Y}-\CJ^{n,k}_{C,Y}),$$
   where the source is an open subscheme of $P^{2g-2-n,\delta-k-1}_\vf$ and
 where the target is the scheme-theoretic image.  In \cite{R80}, Rego
proved versions of some of the preliminary results in Section 4.

Concerning $\pi_\vf$, our main result is Theorem~(6.3), which asserts
that, if $\vf\: C^\dagger\to C$ is the blow-up at a double point,
then $\pi_\vf$ makes $P^{n,1}_\vf$ a (Zariski locally trivial)
$\IP^1$-bundle over $\CJ^n_{\dr C}$.  Theorem~(6.3) extends \cite{AK90,
(10)}, again relaxing the condition that $Y$ is an ordinary node or
cusp.

Throughout,
 we'll work with {\it projective\/} curves and with {\it
projective families\/} of curves because our goal is to further the
theory of the compactified Jacobian.  However, in a number of places,
this hypothesis is unnecessary, provided, when
 needed,
 we assume
that the birational maps at hand are proper.  In this way, with little
or no effort, we can generalize the following results: (3.4), (3.7),
(3.8), (4.2), (5.1), (6.1), (6.2) and (6.4).  To do so, though, we must
also make minor  changes in the preliminary discussions in
(2.1), (2.2), (3.1), and (3.2).

 In outline,
 here's how the preceding results are used in \cite{EGK99} to
prove the autoduality theorem stated above.  First of all, it
is relatively easy, using the determinant of cohomology, to construct a
right inverse to $A_{\cI,C}^*$.  The hard part is to prove that
$A_{\cI,C}^*$ is a monomorphism.  So let $\cN$ be an invertible sheaf on
$\CJ_C^0$.  Suppose that $\cN$ corresponds to a point of
$\Pic_{\CJ_C^0}^0$, and that $A_{\cI,C}^*\cN$ is trivial.  We have to
prove that $\cN$ itself is trivial.

In \cite{EGK99}, we'll prove a form of the theorem of the cube, and use
it to prove that $A_{\cI,C}^*\cN$ is independent of $\cI$.  In fact,
$A_C^*\cN$ is, on $J_C^1\x C$, the pullback of a sheaf $\cN_1$ on
$J_C^1$.  The next step is to prove that $\cN_1$ is trivial.

If $C$ is smooth, then the autoduality theorem is already known.  So
assume that $C$ has a singular point $Q$, which, by hypothesis, is a
double point.  Let $\vf\:\dr C\to C$ be the blowup at $Q$.  Then
$\pi_\vf$ makes $P^{0,1}_\vf$ a $\IP^1$-bundle over $\CJ^0_{\dr C}$ by
Theorem~(6.3).  On each $\IP^1$, the restriction of
 $\kappa_\vf^*\cN$
 is
trivial because $\cN$ corresponds to a point of $\Pic^0_{\CJ_C^0}$.
Hence,
 $\kappa_\vf^*\cN$
 is the pullback of a sheaf $\dr\cN$ on $\CJ^0_{\dr
C}$.

Set $\dr\cI:=\vf^*\cI$.  Using the commutativity of the product square
above and the commutativity of a related square that is also treated in
Corollary (5.5), we can see that $A_{\dr\cI,\dr C}^*\dr\cN$ is equal to
$\vf^*A_{\cI,C}^*\cN$, which is trivial by assumption.  Now, $\dr C$
also has at worst double points by Part (1) of Lemma (6.4).  So by
induction on the genus drop to the normalization, we may assume that
autoduality holds for $\dr C$.  Hence $\dr\cN$ is trivial.  Since
$\pi_\vf^*\dr\cN$ is equal to $\kappa_\vf^*\cN$, the latter is trivial.  So
thanks again to the commutativity of the product square above,
$(1\x\vf )^*A_C^*\cN$ is trivial.  Since $(1\x\vf )^*A_C^*\cN$
 is, on $J^1_C\x\dr C$,
 equal to the pullback of $\cN_1$, the latter is trivial as desired.

Once again, consider the product square above, and recall that $A_C$
is smooth by Corollary (2.6).  Consider the descent data on $(1\x\vf
)^*A_C^*\cN$ with respect to $\Lambda_\vf$;
 since $\kappa_\vf^*\cN$ is trivial,
 this data is trivial.
Hence, so is that on $A_C^*\cN$ with respect to $A_C$.  Thus
$\cN|V$ is trivial, where $V:=A_C(J^1_C\x C)$.
 Since  $\dim(\CJ^n_C-V)\le g-2$
  by Corollary (6.8) and since $\CJ_C^0$ is known to be
of dimension $g$ and
 Cohen--Macaulay
\cite{AIK76, (9), \p.8}, therefore $\cN$ is trivial as desired.  Thus
the autoduality theorem is proved.

In short, in Section~2, we introduce the bigraded Abel maps $\alpha_\cF$
and $A_\cF$, and prove that they are smooth and surjective when
$\cF=\bigomega$ and when $\cF=\cO_C$ and $C$ is Gorenstein,
respectively.  In Section~3, we introduce the presentation schemes
$P^{n,r}_\vf$ and the maps $\kappa_\vf\:P^{n,r}_\vf \to
\CJ^{n+\delta-r}_C$ and $\pi_\vf\:P^{n,r}_\vf \to \CJ^n_{C^\dagger}$,
and we prove some general properties of theirs.  In Section~4, we use
certain Abel maps and presentation schemes to produce two complementary
filtrations of $\CJ^n_C$.  In Section~5, we construct a lifting of the
Abel map $\alpha_\cF$ of bidegree $(m,1)$, and prove that the lifting
completes a product square generalizing the one above.  Finally, in
Section~6, we consider $C$ with double points, and tie together the
material in the preceding sections.  In particular, we show
that, if $\vf\: C^\dagger\to C$ is the blow-up at a double point, then
$\pi_\vf$ makes $P^{n,1}_\vf$ a $\IP^1$-bundle over $\CJ^n_{\dr C}$.  We
also characterize the $C$ having at worst double points as those
Gorenstein $C$ whose Abel map $A_C\:J^1_C\x C\to\CJ^0_C$ has, in its
image, every point of dimension $g-1$, where $g$ is the arithmetic
genus.

\sct2 Abel maps

\dsc1 The compactified Jacobian.  Recall some basic theory from
\cite{AK76}, \cite{AK79II}, or \cite{AK80}.  Let $C/S$ be a {\it flat
projective family of integral curves\/}; by such, we mean that the
scheme $S$ is locally Noetherian, and the structure map $p\:C\to S$ is
flat and projective, and has geometrically integral fibers of dimension
1.  By a (relative) {\it torsion-free rank-$1$} sheaf $\cI$ on $C/S$, we
mean an $S$-flat coherent $\cO_C$-module $\cI$ such that, for each point
$s$ of $S$, the fiber $\cI(s)$ is a torsion-free rank-1 sheaf on the
fiber $C(s)$.  Moreover, we say that $\cI$ is of {\it degree} $n$ if
each $\cI(s)$ satisfies the relation,
        $$\chi(\cI(s))-\chi(\cO_{C(s)})=n;\eqno\Cs1.1)$$
 here $n$ is an integer-valued function that is constant on each
connected component of $S$.

Given $n$, form the \'etale sheaf associated to the ``presheaf,'' or
functor, that assigns to each locally Noetherian $S$-scheme $T$ the set
of isomorphism classes of torsion-free rank-1 sheaves of degree $n$ on
$C_T/T$ (where $C_T$ stands for $C\x T$).  This sheaf is representable
by a projective $S$-scheme, denoted $\CJ^n_{C/S}$, or $\CJ^n_{C}$, or
simply $\CJ^n$.  It (obviously) contains $J^n:=\Pic^n_{C/S}$ as an open
subscheme.

Forming the scheme $\CJ^n$ commutes with changing the base $S$; that is,
for every $S$-scheme $T$, we have a natural isomorphism,
        $$\CJ^n_{C_T/T}=\CJ^n_{C/S}\x T.$$
 Furthermore, if the smooth locus of $C/S$ admits a section, then every
$T$-point of $\CJ^n$ arises from some torsion-free rank-1 sheaf on
$C_T/T$.  In particular, then $C\x\CJ^n$ carries a universal (or
Poincar\'e) sheaf.

 \dsc2 The dualizing sheaf.
 Let $C/S$  be a flat projective family of integral curves.  We'll
denote its (relative) dualizing sheaf by $\bigomega_{C/S}$, or
$\bigomega_C$, or simply $\bigomega$.  Recall some basic facts about
$\bigomega$; see \cite{AK80, (6.5), \p.95} for example. \par
First, forming $\bigomega$
commutes with changing the base; that is, for every $S$-scheme $T$, the
pullback $\bigomega_T$ of $\bigomega$ to $C_T$ is the dualizing sheaf of
$C_T/T$.
 \par
  Second, $\bigomega$ is a torsion-free rank-1 sheaf on $C/S$,
and it is of degree $2g-2$, where $g$ denotes the (locally constant)
arithmetic genus (of the fibers) of $C/S$.

Third, let $T$ be an $S$-scheme, and $\cI$ a torsion-free rank-1 sheaf
on $C_T/T$.  Then the sheaf $\cExt^1_{C_T}(\cI,\bigomega_T)$
vanishes on the fibers.  Therefore, the sheaf
        $$\cI^\star:=\cHom_{C_T}(\cI,\bigomega_T)$$
 is flat over $T$, and forming it commutes with changing $T$ by
\cite{AK80, Thm.~1.10, \p.61}.  Hence, $\cI^\star$ is torsion free and
of rank 1 on $C_T/T$.  Say $\cI$ is of degree $m$.  Then $\cI^\star$
is of degree $2g-2-m$ because, on each fiber $X(t)$, duality yields the
relation,
        $$\chi(\cI(t)^\star)=-\chi(\cI(t)).$$

Fourth, the natural map  is an isomorphism,
        $$\cI\risom(\cI^\star)^\star.\eqno\Cs2.1)$$
Indeed, it suffices to check this statement on each fiber; so we may
assume that $S$ is the spectrum of a field.  Then the map is well known
to be an isomorphism locally at each point of $C$; see \cite{E95,
Thm.\uts 21.21c, p.\uts 538} for example.  So, in particular,
$\bigomega^\star=\cO_C$ since $\cO_C^\star=\bigomega$.

Fifth, if $C$ is Gorenstein, then $\bigomega$ is invertible.  Hence, if
$\cI$ is invertible at $Q\in C_T$, then so is $\cI^\star$.  The converse
follows because of \Cs2.1).

Sending $\cI$ to $\cI^\star$ defines a bijection on $T$-points, so an
isomorphism,
       $$\iota^m\:\CJ^m\risom\CJ^{2g-2-m}.$$
 Its inverse is $\iota^{2g-2-m}$.  So $\iota^m$ is the $m$th component
of a canonical involution of the whole compactified Jacobian
$\coprod \CJ^n$.  We will often
 abbreviate $\iota^m$ by $\iota$.

 \dsc3 Abel maps and tensor products.  Let $C/S$ be a flat projective
family of integral curves.  Let $\cF$ be a torsion-free rank-1 sheaf of
degree $f$ on $C/S$, and $n\ge 0$.  Recall that there is a natural map,
        $$\cA^n_\cF\:\,\Quot^n_\cF\to \CJ^{f-n},$$
 called (the $n$th component of) the {\it Abel map\/} in \cite{AK80,
(8.2), \p.101}; it is defined on $T$-points by sending a short exact
sequence,
        $$0\to\cJ\to\cF_T\to\cG\to0,\hbox{ to }\cJ.$$
We will often abbreviate $\cA^n_\cF$ by $\cA_\cF$.

Particularly important is the case where $\cF=\bigomega$.  In this case,
$\cF^\star=\cO_C$.  So the injection $\cJ\into\cF_T$ corresponds to a
map,
         $$u\:\cO_{C_T}\to\cJ^\star.$$
 On each fiber of $C_T/T$, the support of $\cG$ is finite; hence, on
each fiber, $u$ is nonzero, so injective.  Therefore, $u$ is injective
and its cokernel is flat by \cite{EGA, IV$_3$-11.3.7, \p.135}.
Moreover, if $\cJ^\star$ is invertible, then $u$ corresponds, in turn,
to a relative effective Cartier divisor $D$ whose associated sheaf is
$\cJ^\star$.  Thus the composition,
   $$\iota\circ\cA_\bigomega\:\,\Quot^n_\bigomega\to \CJ^{n},$$
 may be viewed as a natural extension of the traditional Abel map for
the Jacobian, which carries a relative effective Cartier divisor $D$ to
its associated sheaf $\cO(D)$.

Let $m$ and $n$ be arbitrary integers, and let $U^{m,n}$ be the open
subscheme of the product $\CJ^m\x\CJ^n$ that represents the \'etale
sheaf associated to the subfunctor whose $T$-points are the pairs of
torsion-free rank-1 sheaves $(\cI,\,\cJ)$ on $C_T/T$ such that $\cI$
is invertible where $\cJ$ is not.  Then, sending $(\cI,\,\cJ)$ to
$\cI\ox\cJ$ defines a map,
        $$\mu^{m,n}\:U^{m,n}\to\CJ^{m+n};$$
 indeed,  $\cI\ox\cJ$ is a torsion-free rank-1 sheaf because locally
$\cI\ox\cJ$ is isomorphic to $\cI$ or to $\cJ$.  We will often
 abbreviate $\mu^{m,n}$ by $\mu$.

 Assume $n\ge 0$ again.  Let $\cF$ be arbitrary, and form the map and
preimage,
        $$\displaylines{\iota\x(\iota\circ\cA_\cF)\:\,
        \CJ^m\x\Quot^n_\cF\to \CJ^{2g-2-m}\x\CJ^{2g-2-f+n},\cr
        W^{m,n}_\cF:=\bigl(\iota\x(\iota\circ\cA_\cF)\bigr)^{-1}
        U^{2g-2-m,2g-2-f+n}.\cr}$$
 Finally, restrict $\iota\x(\iota\circ\cA_\cF)$ to $W^{m,n}_\cF$, and
follow with $\iota\circ\mu$,  obtaining a map,
        $$\alpha^{m,n}_\cF\:W^{m,n}_\cF \to \CJ^{f-2g+2+m-n}.$$
  We'll call  $\alpha^{m,n}_\cF$  the {\it bigraded Abel} map
 {\it associated to} $\cF$, and often abbreviate $\alpha^{m,n}_\cF$
by $\alpha_\cF$.

Note that, if a $T$-point $t$ of $W^{m,n}_\cF$ represents a pair of
sheaves $\cI$ and $\cJ$ on $C_T/T$, where $\cJ$ is a subsheaf of the
pullback $\cF_T$, then $\cI^\star$ is invertible where $\cJ^\star$ is
not, and $\alpha_\cF(t)$ represents the sheaf $\cK$ defined by the
formula,
        $$\cK:=({\cI}^\star\ox{\cJ}^\star)^\star.\eqno\Cs3.1)$$

 Another important case occurs when $\cF=\cO_C$.  Then it is customary
to write $\Hilb^n_C$ instead of $\Quot^n_{\cO_C}$.  Let's also
abbreviate
$W^{m,n}_{\cO_C}$
 by $W^{m,n}_C$, and
$\alpha^{m,n}_{\cO_C}$
by $\alpha^{m,n}_C$ or $\alpha_C$.

 \thm4 For any flat projective family of integral curves, the bigraded
Abel map associated to $\bigomega$,
     $$\alpha_\bigomega\:\,W^{m,n}_\bigomega \to \CJ^{m-n},$$
 is smooth and surjective.
 \pf Let $C/S$ be the family.  To prove that $\alpha_\bigomega$ is
surjective, we may assume that $S$ is the spectrum of an algebraically
closed field.  Now, given a closed point of $\CJ^{m-n}$, let $\cK$ be a
corresponding sheaf on $C$.  Let $Q$ be a simple closed point of $C$,
and set $\cJ:=\bigomega(-nQ)$ and $\cI:=\cK(nQ)$.  Then
  $${\cI}^\star\ox{\cJ}^\star=\cK^\star(-nQ)\ox\cO_C(nQ)=\cK^\star.$$
 Thus $\alpha_\bigomega$ is surjective.

Let's now prove
 that $\alpha_\bigomega$ is smooth.  It is if it is
 after we change the base via a smooth surjective map $S'\to S$.  So we
may assume that the smooth locus of $C/S$ has a section.

Let's use the infinitesimal criterion for smoothness
\cite{EGA, IV$_4$-17.5.4, \p.69}.
 So
 assume that $S$ is the spectrum of an Artin local
ring.  Let $S'$ be a closed subscheme of $S$ whose ideal is of length 1.
Set $C':=C\x S'$ and $\bigomega':=\bigomega|C'$.
Given an $S'$-point $w'$ of $W^{m,n}_\bigomega$
whose image under $\alpha_\bigomega$ lifts to an $S$-point $k$ of
$\CJ^{m-n}$, we must construct an $S$-point $w$ of $W^{m,n}_\bigomega$
such that $w|S'=w'$ and $\alpha_\bigomega(w)=k$.

Since the smooth locus of $C/S$ has a section, $w'$ arises from a pair
of sheaves $\cI',\ \cJ'$ on $C'/S'$, where $\cJ'$ is a subsheaf of
$\bigomega'$ with flat quotient.  Similarly, $k$ arises from a sheaf
$\cK$ on $C/S$.  In these terms, we have to lift $\cI'$ to a flat sheaf
$\cI$ on $C/S$ and to lift $\cJ'$ to a subsheaf $\cJ$ of $\bigomega$
such that Formula
 \Cs3.1)
 holds.  Note that $\bigomega/\cJ$ is
automatically flat by \cite{EGA, IV$_3$-11.3.7, \p.135}.

Since $C'$ is of dimension one, ${\cJ'}^\star$ is invertible off a
finite set $Y_1$.  But ${\cI'}^\star$ is invertible at each point of
$Y_1$ since the pair $(\cI',\,\cJ')$ corresponds to the $S$-point $w'$ of
$W^{m,n}_\bigomega$.  So ${\cI'}^\star$ is trivial on some affine open
subscheme $U'_1$ of $C'$ containing $Y_1$; fix a trivialization
$\xi'_1\:\cO_C|U'_1\risom{\cI'}^\star|U'_1$.  The complement $Y_2$ of
$U'_1$ is also finite, and ${\cJ'}^\star$ is invertible at each point of
$Y_2$, because $U'_1$ contains $Y_1$.  Hence ${\cJ'}^\star$ is trivial
on some affine open subscheme $U'_2$ containing $Y_2$; fix a trivialization
$\xi'_2\:\cO_C|U'_2\risom{\cJ'}^\star|U'_2$.

Let's first lift $\cI'$.  Tensoring $\xi'_2$ by ${\cI'}^\star|U'_2$ yields
an isomorphism
$\eta'_2\:{\cI'}^\star|U'_2\risom\cK^\star|U'_2$.  Set
$U'_{12}:=U'_1\cap U'_2$.  Define
$\xi'_{12}\:\cO_C|U'_{12}\to\cK^\star|U'_{12}$ by
$\xi'_{12}:=\eta'_2\xi'_1$.  Let $U_1$ and $U_2$ denote the affine open
subschemes of $C$ on the same underlying sets as $U'_1$ and $U'_2$.
Set $U_{12}:=U_1\cap U_2$. Since $\cK^\star$ is flat and since $U_{12}$
is affine, $\xi'_{12}$ lifts to a trivialization
$\xi_{12}$ of $\cK^\star|U_{12}$.
Use $\xi_{12}$ to patch together
$\cO_{C}|U_1$ and $\cK^\star|U_2$, and then apply the operator
$\star$.  The result is a lifting $\cI$ of $\cI'$.  In addition,
we obtain isomorphisms $\xi_1\:\cO_C|U_1\risom\cI^\star|U_1$ and
$\eta_2\:\cI^\star|U_2\risom\cK^\star|U_2$
lifting $\xi'_1$ and $\eta'_2$, respectively.

Of course,  $\cI$ is not unique.  However, the difference between it and
any other lifting is represented by a vector in
 $${\rm H}^1\bigl(C'', \cHom_{C''}({\cI''}^\star,{\cI''}^\star)\bigr);$$
 here and below, ${}''$ means the restriction over the reduced point of $S$.
We will find a suitable vector to use to modify $\cI$ so that we
can construct a suitable lifting $\cJ$ of $\cJ'$.

The inclusion of $\cJ'$ in $\bigomega'$ corresponds to an injection
$u'\:\cO_{C'}\into{\cJ'}^\star$, whose cokernel is flat; see \Cs3). Let
$a'\in(U'_2,\cO_{C'})$ be such that $u'|U'_2=a'\xi'_2$.  Since
$U_2$ is affine, $a'$ lifts to an element $a\in\Gamma(U_2,\cO_C)$.  Use $a$
to define a divisor $D$ of $C$ with support in $U_2$; note that $a$
is regular and that $D$ is flat by \cite{EGA, IV$_3$-11.3.7, \p.135},
because $D|U'_2$ is just the divisor associated to the injection
$u'|U'_2$.  Set $D':=D|C'$.

Set $v':=u'\ox 1_{{\cI'}^\star}$; so $v'\:{\cI'}^\star\to{\cK'}^\star$.
By construction, $v'|U'_2=a'\eta'_2$; so $v'|U'_2$ is injective. Since
${\cI'}^\star|U'_1$ is trivial, $v'|U'_1$ is injective with flat
cokernel, because $u'|U'_1$ is so. Thus $v'$ factors through an
injection,
        $$v'_0\:{\cI'}^\star\to{\cK'}^\star(-D'),$$
 whose cokernel is flat,
because the cokernel of $v'|U'_1$ is so, and because
$v'_0|U'_2$ is bijective.

Consider the obstruction to lifting $v'_0$ to a map from $\cI^\star$ to
$\cK^\star(-D)$. This obstruction is represented by a vector,
 $$e\in\Ext^1_{C''}\bigl({\cI''}^\star,{\cK''}^\star(-D'')\bigr).$$
 The ``local obstruction'' is the image of $e$ in
  $${\rm H}^0\bigl(C'', \cExt^1_{C''}(
        {\cI''}^\star,{\cK''}^\star(-D''))\bigr),$$
 and it vanishes; indeed, it vanishes on $U_1''$ because $\cI''$ is
trivial there, and it vanishes on $U_2''$ because
$v'|U'_2$ lifts to $a\eta_2$.  Hence,
$$e\in {\rm H}^1\bigl(C'',
         \cHom_{C''}({\cI''}^\star,{\cK''}^\star(-D''))\bigr)
 \subset \Ext^1_{C''}\bigl({\cI''}^\star,{\cK''}^\star(-D'')\bigr).$$

Since $v'_0$ is injective and its cokernel is flat, its restriction to
$C''$ is injective.  So the restriction's cokernel has finite support.
Hence the restriction induces a surjection,
  $${\rm H}^1\bigl(C'',\cHom_{C''}({\cI''}^\star,{\cI''}^\star)\bigr)
 \onto{\rm H}^1\bigl(C'',
  \cHom_{C''}({\cI''}^\star,{\cK''}^\star(-D''))\bigr).$$
 Lift $e$ and modify $\cI$.  Then $e$ vanishes, so
$v'_0$  lifts to a map $v_0\:\cI^\star\to\cK^\star(-D)$.  Composing
$v_0$ with the inclusion of $\cK^\star(-D)$ in $\cK^\star$, we obtain a
map $v\:\cI^\star\to\cK^\star$ lifting $v'$.

Finally, define a subsheaf $\cJ$ of $\bigomega$ as follows.  On $U_1$,
define $\cJ$ to be the image of the injection
$\cK|U_1\into\bigomega|U_1$ corresponding under the operator $\star$ to
the composition $v\xi_1$.  On $U_2$, define $\cJ$ to be $\bigomega(-D)$.
Then $\cJ$ lifts $\cJ'$; it does so on $U_1$ because $\cI^\star$ is
trivial there, and it does so on $U_2$ because
$\cJ'|U'_2=\bigomega'(-D')|U'_2$
by construction of $D'$.
 Moreover, $v\:\cI^\star\to\cK^\star$ factors as follows:
        $$v\:\cI^\star\to\cI^\star\ox\cJ^\star\risom\cK^\star,$$
where the first map is obtained by tensoring the injection
$\cO_C\into\cJ^\star$ with $\cI^\star$, and the second is an
 induced isomorphism.
The theorem is now proved.

 \cor5 Let $C/S$ be a flat projective family of integral curves, and
assume that $C/S$ is Gorenstein.  Then the bigraded Abel map associated
to $\cO_C$,
     $$\alpha_C\:W^{m,n}_C \to \CJ^{m-n-2g+2},$$
 is smooth and surjective.
 \pf Since
 $\bigomega$ is invertible, tensoring with it defines two isomorphisms,
        $$\tau\:\CJ^m\risom\CJ^{m+2g-2} \and
         \sigma\:\Hilb^n_C\risom\Quot^n_\bigomega.\eqno\Cs5.1)$$
 Moreover, an arbitrary sheaf $\cL$ on $C$ is invertible at a given
point if and only if $\cL\ox\bigomega$ is so.  It follows that
        $$(\tau\x\sigma)(W^{m,n}_C)=W^{m+2g-2,n}_\bigomega \and
        \alpha_\bigomega\circ(\tau\x\sigma)
        =\tau\circ\tau\circ\alpha_C.$$
 The assertion now results from Theorem~\Cs4).

 \cor6 Let $C/S$ be a flat projective family of integral curves, and
assume that $C/S$ is Gorenstein.  Then the Abel
map,
        $$A_C\:W^{m,n}_C \to \CJ^{m-n}, \hbox{ defined by }
        (\cI,\cJ)\mapsto\cI\ox\cJ,$$
   is well-defined, smooth, and surjective.
 \pf
 Thanks to Corollary \Cs5), we need only prove that, with $\tau$ as in
\Cs5.1),
        $$A_C=\tau\circ\alpha_C.\eqno\Cs6.1)$$

First, let's make a simple general observation: given any three
quasi-coherent sheaves $\cM$, $\cN$, and $\cP$ on any scheme, the
canonical map is an isomorphism,
    $$\cHom(\cM,\cP)\ox\cN\risom\cHom(\cM,\cP\ox\cN),\eqno\Cs6.2)$$
 if $\cM$ is invertible where $\cN$ isn't.  Indeed, the question is
local, and the map is functorial in $\cM$ and $\cN$.  However, the
question is trivial if either $\cM$ or $\cN$ is trivial.

Whether $C/S$ is Gorenstein or not, adjoint associativity and \Cs2.1)
yield
        $$({\cI}^\star\ox{\cJ}^\star)^\star
        =\cHom({\cI}^\star,\cHom({\cJ}^\star,\bigomega_T))
        =\cHom({\cI}^\star,\cJ).\eqno\Cs6.3)$$
 Now, $(\cI,\cJ)$ defines a $T$-point of $W^{m,n}_C$; so $\cI^\star$ is
invertible where $\cJ^\star$ is not.  Since $C/S$ is Gorenstein,
$\cJ^\star$ is invertible precisely where $\cJ$ is, and $\bigomega$ is
invertible everywhere; see \Cs2).  Hence, \Cs6.2)
 and \Cs2.1) yield
        $$\cHom({\cI}^\star,\cJ)
        =\cHom({\cI}^\star,\bigomega_T)\ox(\cJ\ox\bigomega_T^{-1})
        =\cI\ox\cJ\ox\bigomega_T^{-1}.\eqno\Cs6.4)$$
 Thus \Cs6.1) holds, and the corollary is proved.

\sct3 Presentation schemes

 \dsc1 Existence.  (See \cite{AK90}.)  Let $p\:C\to S$ and $\dr p\:\dr
C\to S$ be flat projective families of integral curves, and $\vf\:\dr
C\to C$ a (relative) {\it birational\/} $S$-map; that is, $\vf$ is an
$S$-map whose fibers $\vf(s)\:\dr C(s)\to C(s)$ are all birational.
Then $\vf$ is proper and quasi-finite, so finite.  Its comorphism
$\vf^\sharp\:\cO_C\to \vf_*\cO_{\dr C}$ has injective fibers; so
$\vf^\sharp$ is injective and its cokernel is flat.  Set
        $$\delta_\vf:= \length\bigl((\vf_*\cO_{\dr C}/\cO_C)(s)\bigr).$$
 By flatness, $\delta_\vf$ is locally constant.  We'll call $\delta_\vf$
 the {\it genus drop\/} from $C$ to $\dr C$ or along $\vf$, and we'll
often abbreviate $\delta_\vf$ by $\delta$.

In $C$ and $\dr C$, form the sets,
        $$Y:=\Supp\bigl(\vf_*\cO_{\dr C}/\cO_C\bigr)\and
          \dr Y:=\vf^{-1}Y,$$
 and give them the subscheme structures defined by the annihilator
$\cI_Y$ of $\vf_*\cO_{\dr C}/\cO_C$.  We'll call $Y\subset C$ and $\dr
Y\subset \dr C$ the {\it conductor subschemes}, and $\cI_Y$ the {\it
conductor ideal} of $\dr C/C$ or $\vf$.  Note that
 evaluating at 1 induces an isomorphism,
        $$\cHom(\vf_*\cO_{\dr C},\cO_C)\risom\cI_Y.$$
 Moreover,
 if $T$ is an
$S$-scheme, then $Y_T$ contains the conductor subscheme in $C_T$ of
 $\vf_T$,
 and these two subschemes of $C_T$ have the same underlying set.

If $\delta=1$, then $p$ restricts to an isomorphism $Y\risom S$, and
$\dr Y$ is $S$-flat of relative length 2.  Furthermore, forming $Y$ and
$\dr Y$ commutes with changing $S$.
  (See \cite{AK90, Prop.~2, \p.17}.)

By duality theory, $(\vf_*\cO_{\dr C})^\star=\vf_*(\cO_{\dr C}^\star)$.
Since $\bigomega_C=\cO_C^\star$ and $\bigomega_{\dr C}=\cO_{\dr C}^\star$,
the comorphism $\vf^\sharp$ gives rise to a map,
        $$c_\vf\: \vf_*\bigomega_{\dr C}\to\bigomega_C.$$
Furthermore, $c_\vf$ is injective, $\cCok(c_\vf)$ is a flat $\cO_Y$-module
of relative length $\delta$, and forming $c_\vf$ commutes with changing $S$.

Suppose $C/S$ is Gorenstein along $Y$.  Then $\bigomega_C$ is invertible
wherever $\vf_*\cO_{\dr C}$ isn't.  So
        $$\cHom(\vf_*\cO_{\dr C},\bigomega_C)=
        \cHom(\vf_*\cO_{\dr C},\cO_C)\ox\bigomega_C$$
 by the general fact explained in the proof of Corollary~(2.6).  Hence
  $$\vf_*\bigomega_{\dr C}=\vf_*(\cO_{\dr C}^\star)=(\vf_*\cO_{\dr C})^\star
        =\cHom(\vf_*\cO_{\dr C},\cO_C)\ox\bigomega_C
        =\cI_Y\ox\bigomega_C,\eqno\Cs1.1)$$
 and so $\cCok(c_\vf)=\cO_Y\ox\bigomega_C$.
Therefore, $Y$ is flat of relative
length $\delta$, and forming $Y$ commutes with changing $S$.  Also,
$\dr Y$ is flat of relative length $2\delta$, and forming $\dr Y$
commutes with changing $S$ because of the short exact sequence,
        $$0\to\cO_C/\cI_Y\to \vf_*\cO_{\dr C}/\cI_Y\to
        \vf_*\cO_{\dr C}/\cO_C\to0.$$

Let $T$ be an $S$-scheme, and $\dr\cI$ a torsion-free rank-1 sheaf on
$\dr C_T/T$.  By definition, a {\it presentation\/} of $\dr\cI$ for
$\vf$ is an injective map,
                $$h\:\cI \to \vf_{T*}\dr\cI,$$
 (where $\vf_T:=\vf\x1_T$) such that $\cCok(h)$ is a $T$-flat
$\cO_{Y_T}$-module.

If $\cCok(h)$ has relative length $r$ over $T$ (resp., if $\dr\cI$ is of
degree $n$), then we'll say that $h$ is {\it of colength} $r$ (resp.,
{\it of degree} $n$).  If so, then $\cI$ is a torsion-free rank-1 sheaf
on $C_T/T$, and it is of degree $n+\delta-r$.  For example, the
comorphism $\vf^\sharp\:\cO_C\to \vf_*\cO_{\dr C}$ is a presentation of
colength $\delta$.

Since $Y_T$ contains the conductor subscheme of $\vf_T$, any presentation
for $\vf_T$ is a presentation for $\vf$.  The converse holds if $\delta=1$
or if $C/S$ is Gorenstein along $Y$, because then forming $Y$ commutes
with changing $S$, as noted above.

Two presentations, $h_1\:\cI_1 \to \vf_{T*}\dr\cI_1$ and $h_2\:\cI_2 \to
\vf_{T*}\dr\cI_2$, are considered {\it equivalent\/} if there exist an
invertible sheaf $\cN$ on $T$ and a commutative diagram
 with vertical isomorphisms,
 $$ \CD
    \cI_1 @>h_1>> \vf_{T*}\dr\cI_1\\
     @VV\simeq V              @VV\simeq V\\
 \cI_2\ox p_T^*\cN @>h_2\ox p_T^*\cN>> \vf_{T*}\dr\cI_2\ox p_T^*\cN
 \endCD $$

For each $n$ and $r$, there exists a projective $S$-scheme $P_\vf^{n,r}$
that parameterizes the presentations of colength $r$ of sheaves $\dr
\cI$ of degree $n$ up to equivalence; the existence of $P_\vf^{n,r}$ is
proved in two different ways in \cite{AK90, Thm.~12, \p.22} for
$r=\delta$, but both proofs work for any $r$.  More precisely, if there
exists a universal sheaf on $\dr C\x\CJ^n_{\dr C}$, where
$\CJ^n_{\dr C}$
is the compactified Jacobian, then the $T$-points of $P_\vf^{n,r}$
(or $S$-maps $T\to P_\vf^{n,r}$) correspond bijectively to the equivalence
classes of presentations of colength $r$ and degree $n$ for $\vf_T$; in
other words, then $P_\vf^{n,r}$ represents the functor of equivalence
classes of presentations.  In general, $P_\vf^{n,r}$ represents the
associated \'etale sheaf.

Forming $P_\vf^{n,r}$ commutes with changing $S$ if $\delta=1$ or if
$C/S$ is Gorenstein along $Y$, because, as noted above, then any
presentation for $\vf_T$ is a presentation for $\vf$, and conversely.

 \dsc2 Copresentations.
 Let $\vf\:\dr C\to C$ be a birational $S$-map of flat projective
families of integral curves.
 Let $T$ be an $S$-scheme, and $\dr\cI$ a torsion-free rank-1 sheaf on
$\dr C_T/T$.  By a {\it copresentation\/} of $\dr\cI$ for $\vf$, we'll
mean an injective map,
                $$c\:\vf_{T*}\dr\cI \to\cI,$$
 such that $\cI$ is a torsion-free rank-1 sheaf on $C_T/T$, and
$\cCok(c)$ is a $T$-flat $\cO_{Y_T}$-module. If $\cCok(c)$ has relative
length $r$ (resp., if $\dr\cI$ is of degree $n$), then we'll say that
$c$ is {\it of colength} $r$ (resp., {\it of degree} $n$).  For example,
the natural map $c_\vf\:\vf_*\bigomega_{\dr C}\to\bigomega_C$ is a
copresentation of colength $\delta$.

The operator $\star$ of (2.2) commutes with $\vf_{T*}$; indeed, duality
theory provides a canonical isomorphism,
        $$(\vf_{T*}{\dr\cI})^\star=\vf_{T*}({\dr\cI}^\star).\eqno\Cs2.1)$$
 So  $c$ gives rise to a map,
        $$c^\star\:\cI^\star\to \vf_{T*}({\dr\cI}^\star).$$
 Its source and target are torsion-free rank-1 sheaves on $ C_T/T$,
and forming $c^\star$ commutes with changing $T$.  Hence, $c^\star$ is
injective, and $\cCok(c^\star)$ is $T$-flat.

The cokernels of $c$ and $c^\star$ are also related by duality theory:
        $$\cCok(c^\star)=\cExt^1(\cCok(c),\bigomega_{C_T})
        \and \cCok(c)=\cExt^1(\cCok(c^\star),\bigomega_{C_T}).$$
 Hence the two cokernels have the same annihilator.  So
$\cCok(c^\star)$ is an $\cO_{Y_T}$-module.  Thus $c^\star$ is a
presentation of ${\dr\cI}^\star$ of the same colength as $c$.  For
example, $c_\vf^\star$ is equal to the comorphism $f^\sharp\:\cO_C\to
f_*\cO_{\dr C}$.

Similarly, every presentation $h\:\cI \to \vf_{T*}\dr\cI$ gives rise to a
map,
        $$h^\star\: \vf_{T*}({\dr\cI}^\star)\to\cI^\star,$$
 which is a copresentation of the same colength.  Moreover,
$h^{\star\star}=h$ and $c^{\star\star}=c$.  Furthermore, the notion of
equivalence of presentations gives rise to a notion of equivalence of
copresentations.  Therefore, for each $n$ and $r$, the presentation
scheme $P_\vf^{n,r}$ parameterizes the copresentations of colength $r$ of
sheaves $\dr\cI$ of degree $2g-2\delta-2-n$ up to equivalence,
where $g$ denotes the
(locally constant) arithmetic genus of $C/S$.

\dsc3 Maps from  $P_\vf^{n,r}$.
 Let $\vf\:\dr C\to C$ be a birational $S$-map of flat projective
families of integral curves.
There are natural maps from $P_\vf^{n,r}$ to the compactified Jacobians of
$C/S$ and $\dr C/S$:
        $$\kappa^{n,r}_\vf\:P_\vf^{n,r}\to\CJ^{n+\delta-r}_C\and
          \pi^{n,r}_\vf\:P_\vf^{n,r}\to\CJ^n_{\dr C}.$$
 They are defined by sending a presentation $h\: \cI \to \vf_{T*}\dr\cI$
to $\cI$ and to $\dr\cI$, respectively.  They are well defined because
they respect equivalence; $\pi_\vf^{n,r}$ does so thanks to Lemma \Cs4)
below (it is a special case of \cite{AK90, Prop.~3, \p.18}, with a
simpler proof).  We will often
 abbreviate $\kappa^{n,r}_\vf$ and $\pi^{n,r}_\vf$ by $\kappa_\vf$ and
$\pi_\vf$.

 It follows from the definitions that $\pi_\vf^{n,0}$ is an isomorphism.
So the composition,
  $$\epsilon_\vf^n:=\kappa^{n,0}_\vf(\pi^{n,0}_\vf)^{-1}\:\CJ^n_{\dr C} \to
          \CJ^{n+\delta}_C,$$
 is defined; it carries a sheaf $\dr\cI$ on $\dr C_T/T$ to
$\vf_{T*}\dr\cI$.
 We'll often abbreviate $\epsilon_\vf^n$ by $\epsilon_\vf$.

Lemma \Cs4) implies that $\epsilon_\vf$ is a monomorphism.  Since it is
also proper, it is a closed embedding.  So in \cite{AK90, 5, \p.19},
$\epsilon_\vf$ was termed the {\it canonical embedding}.  It commutes
with the canonical involution of (2.2); indeed, Equation \Cs2.1) says
that $\iota\epsilon_\vf=\epsilon_\vf\iota$.

 \lem4
 Let $\vf\:\dr C\to C$ be a birational $S$-map of flat projective
families of integral curves.  Let $\dr\cI_1$ and $\dr\cI_2$ be two
 torsion-free rank-$1$ sheaves on $\dr C/S$, and let $u\:\vf_*\dr\cI_1
\to\vf_*\dr\cI_2$ be a map.  Then there is a unique map
 $\dr u\:\dr\cI_1 \to \dr\cI_2$ such that $\vf_*\dr u=u$.
 \pf Consider the following diagram on $\dr C$:
        $$\CD
        \vf^*\vf_*\dr\cI_1 @>\vf^*u>> \vf^*\vf_*\dr\cI_2 \\
           @VVV                 @VVV\\
           \dr\cI_1     @.      \dr\cI_2
        \endCD$$
  where the vertical maps are the canonical maps.  They are surjective
since $\vf$ is finite.  Hence, since their targets are flat, forming
their kernels commutes with passing to the fibers.  On the fibers, their
kernels have finite support because the fibers of $\vf$ are birational.
 Hence the left kernel maps to zero in $\dr\cI_2$; indeed, the latter is
flat and torsion-free on the fibers of $\dr C/S$, so its associated
points are the generic points of their fibers.
Hence there is a unique map $\dr u\:\dr\cI_1
 \to \dr\cI_2$ making the above diagram commutative.  By adjunction,
$\vf_*\dr u=u$, and the proof is complete.

 \prp5
 Let $\vf\:\dr C\to C$ be a birational $S$-map of flat projective
families of integral curves.  If $\dr C/S$ is smooth along the
conductor scheme $\dr Y$, then,
  locally in the \'etale topology on
$\CJ^n_{\dr C/S}$, the map $\pi_\vf^{n,r}$ is isomorphic to the following
second projection:
        $$\Quot^r_{\vf_*\cO_{\dr Y}/Y/S}\x\CJ^n_{\dr C/S} \to
          \CJ^n_{\dr C/S}.$$
  \pf Set $\dr\CJ:=\CJ^n_{\dr C/S}$.  We may replace $S$ by an \'etale
covering; so we may assume that there is a universal torsion-free
rank-1 sheaf $\dr\cI$ on $\dr C\x\dr\CJ/\dr\CJ$.  By
hypothesis, $\dr Y$ is contained in the smooth locus of $\dr C/S$; so
the restriction $\cF$ of $\dr\cI$ to $\dr Y\x\dr\CJ$ is
invertible.  Let $p_2\:\dr Y\x\dr\CJ\to\dr\CJ$ denote
the projection.  Let $s\in \dr\CJ$.  Since the fiber
$p_2^{-1}(s)$ is finite, there is an open subscheme $U\subseteq\dr
Y\x\dr\CJ$ containing $p_2^{-1}(s)$ such that $\cF|U$ is free.
Set $V:=\dr\CJ-p_2(\dr Y\x\dr\CJ- U)$.  Since $p_2$
is proper, $V$ is an open neighborhood of $s$ in $\dr\CJ$. By
 construction,
 $p_2^{-1}V\subseteq U$.  Hence $\cF|p_2^{-1}V$ is free.

Set $\vf':=\vf\x1_{\dr\CJ}$ and $\cG:=(\vf'_*\dr\cI)|(Y\x\dr\CJ)$.  Over
$\dr\CJ$, the presentation scheme $P_\vf^{n,r}$ is isomorphic to
$\Quot^r_{\cG/Y\x\dr\CJ/\dr\CJ}$; this statement is proved in
\cite{AK90, Prop.~9, \p.20} for $r=\delta$, but the proof works for any
$r$.  Let $\psi\:\dr Y\x\dr\CJ\to Y\x\dr\CJ$ denote the restriction of
$\vf'$.  Since $\cG=\psi_*\cF$, and $\cF|p_2^{-1}V$ is free, there is an
isomorphism,
        $$\cG|(Y\x V)\cong(\psi_*\cO_{\dr Y\x\dr\CJ})|(Y\x V).$$
 Hence $(\pi_\vf^{n,r})^{-1}V$ is $V$-isomorphic to
$\Quot^r_{\vf_*\cO_{\dr Y}/Y/S}\x V$.  The proof is now complete.

\prp6
 Let $\vf\:\dr C\to C$ be a birational $S$-map of flat projective
families of integral curves.
 If $\dr C/S$ is smooth along its conductor subscheme $\dr Y$, then
every $\kappa_\vf^{n,r}$ is finite.
 \pf Since $\kappa_\vf$ is proper, it is finite if it is quasi-finite.
Let $t$ be a geometric point of $\CJ^{n+\delta-r}_C$, representing a
sheaf $\cI$ on the fiber $C(t)$.
  Form $\dr\cI:=\cI\cO_{\dr C(t)}$, the quotient of $\cI\ox\cO_{\dr
C(t)}$ by its torsion subsheaf.
  Consider the natural map $h\:\cI\to \vf(t)_*\dr\cI$.

Let $u\:\cI\to \vf(t)_*\dr\cJ$ be a presentation.  Its adjoint
$\vf(t)^*\cI\to\dr\cJ$ factors through an injection $\dr
u\:\dr\cI\to\dr\cJ$.  So $u$ factors as follows:
        $$u\:\cI@>h>>\vf(t)_*\dr\cI@>\vf(t)_*\dr u>> \vf(t)_*\dr\cJ.$$
 Since $\cCok(u)$ is an $\cO_{Y(t)}$-module, and since $\dr C(t)$ is
smooth along $\dr Y(t)$, there is an effective divisor $\dr D$ on $\dr
C(t)$ such that $\dr u(\dr\cI)=\dr\cJ(-\dr D)$ and such that $\dr
D\le[\dr Y(t)]$, where $[\dr Y(t)]$ is the fundamental cycle of the
scheme $\dr Y(t)$.

 Given another presentation $u_1\:\cI\to \vf(t)_*\dr\cJ_1$, let
$\dr u_1\:\dr\cI\to\dr\cJ_1$ be the induced injection, and $\dr D_1$ the
divisor such that $\dr u_1(\dr\cI)=\dr\cJ_1(-\dr D_1)$. If $\dr D_1=\dr D$,
then there is an isomorphism $v\:\dr\cJ\risom\dr\cJ_1$ such that
$\dr u_1=v\dr u$; hence, $u$ is equivalent to $u_1$.  Now, there are at most
finitely many effective divisors bounded by $[\dr Y(t)]$.  Hence the fiber
of $\kappa_\vf$ over $t$ is finite.  Thus the proposition is
proved.

 \lem7 Let $C$ be a projective integral curve over an algebraically
closed field, $Q\in C$ a singular point, and $\cM\subset\cO_C$ its
ideal.  Form the $\cO_C$-algebra $\cB:=\cHom(\cM,\cM)$ and set
$C^*:=\Spec(\cB)$.
Then these assertions hold.
 \part 1 We have $\cB=\cHom(\cM,\,\cO_C)$.
 \part 2 The structure map  $\psi\:C^*\to C$ is birational with genus drop
$1$ if and only if $C$ is Gorenstein at $Q$.
 \part 3 Assume that $C$ is Gorenstein at $Q$, and let $\vf\:\dr
C\to C$ be a birational map that is nontrivial at $Q$.  Then $\vf$
factors through $\psi$.
 \pf
  Consider (1).  The inclusion $\cM \into \cO_C$ induces an inclusion,
        $$\cB\into\cHom(\cM,\cO_C).$$
 It is an equality.  Indeed, consider a map $u\:\cM_Q\to\cO_{C,Q}$.  If
$u\notin\cB_Q$, then, since $\cM_Q$ is the maximal ideal, $u$ is
surjective, so an isomorphism.  If so, then $\cM_Q$ is invertible.
However, $Q$ is a singular point.  Hence $u\in\cB_Q$.  Thus (1) holds.

Consider (2).  The basic exact sequence $0\to\cM\to\cO_C\to\cO_Q\to0$
yields
        $$0\to\cHom(\cO_C,\cO_C)\to\cHom(\cM,\cO_C)\to
        \cExt^1(\cO_Q,\cO_C)\to0\eqno\Cs7.1)$$
 because $\cHom(\cO_Q,\cO_C)$ and $\cExt^1(\cO_C,\cO_C)$ vanish.  The
first term of \Cs7.1) is equal to $\cO_C$.  The second term is equal to
$\cB$ by (1).  The third term is of length 1 if and only if $C$
is Gorenstein at $Q$.  Hence (2) holds.

Consider (3).  Set $\cA:=\vf_*\cO_{\dr C}$.  The comorphism
$\vf^\sharp\:\cO_C\to\cA$ induces an injection,
        $$u\:\cHom(\cA_Q,\cO_{C,Q})\into\cO_{C,Q}.$$
  If there were a map
$v\:\cA_Q\to\cO_{C,Q}$ such that $v\vf^\sharp_Q(1)
\notin\cM_Q$, then $v$ would be bijective, and so $\vf^\sharp_Q$ would be
bijective too; hence, $u(\cHom(\cA_Q,\cO_{C,Q}))\subset\cM_Q$.  So $u$
induces an injection,
        $$\cHom(\cM_Q,\cO_{C,Q})\into
        \cHom(\cHom(\cA_Q,\cO_{C,Q}),\cO_{C,Q}).$$
 Its source is equal to $\cB_Q$ by (1).  Its target is equal to $\cA_Q$ if
$C$ is Gorenstein at $Q$.  Thus (3) holds.

 \lem8
Let $C$ be a projective integral curve over an algebraically
closed field, and $Q$ a singular point at which $C$ is
Gorenstein.  Let $\cM\subset\cO_C$ denote the ideal of $Q$.  Set
$\cB:=\cHom(\cM,\cM)$ and $C^*:=\Spec(\cB)$.  Let
$\psi\:C^*\to C$ denote the structure map.  Finally,
let $\cI$ be a
torsion-free rank-$1$ sheaf on $C$.  \par
  Then there exists a torsion-free rank-$1$ sheaf $\cI^*$ on $C^*$ such
that $\psi_*\cI^*=\cI$ if and only if $\cI$ is not invertible at $Q$.
If $\cI^*$ exists, then it is unique.
 \pf
 If $\cI^*$ exists, then it is unique by Lemma~\Cs4).  Now, $\cI^*$
exists if and only if $\cI$ is a $\cB$-module, so if and only if
$\cHom(\cI,\cI)$ contains $\cB$, where both $\cO_C$-modules are viewed
as fractional ideals.  Hence, $\cI^*$ does not exist if $\cI$ is
invertible at $Q$, because then $\cHom(\cI,\cI)$ is equal to $\cO_C$ at
$Q$, whereas $\delta_\psi=1$ by Part~(2) of Lemma~\Cs7).  So suppose
$\cI$ is not invertible at $Q$.

For convenience, set $A:=\cO_{C,Q}$ and $M:=\cM_Q$.  Set $B:=\cB_Q$ and
$I:=\cI_Q$.  We have to show that $BI\subset I$.  Set $I':=\Hom_A(I,A)$.
Then $\Hom_A(I',A)=I$ because $A$ is Gorenstein.  Let $\bar A$ denote the
normalization of $A$, and set $K:=\Hom_A(\bar A,A)$.  View $K$ and $I'$
as fractional ideals too.

Since $I'\bar A$ is principal, there is an element $g$ in $I'$ such that
$I'\bar A=g\bar A$.  Then we have $gA\subset I'\subset g\bar A$.  Apply
$\Hom_A(-,A)$ and multiply by $g$.  We get
        $$K\subset gI\subset  A.\eqno\Cs8.1)$$
 Since $I$ is not invertible, $gI\neq A$; hence, $gI\subset M$.

Let $\?M$ denote the radical of $\bar A$.  Then $M\subset \?M$, so
        $$gI\subset \?M.\eqno\Cs8.2)$$
 Since $\bar A$ is normal, $\?M$ is invertible; let $\?M^{-1}$
denote the inverse fractional ideal.  Set $G:=\?M^{-1}K$.
Multiplying \Cs8.2) by $G$ yields $gIG\subset K$.  So \Cs8.1)
yields $gIG\subset gI$.  Therefore,
         $$IG\subset I.\eqno\Cs8.3)$$

Note that $G\not\subset A$.  Indeed, otherwise $G\subset K$ because $K$
is the largest $\bar A$-submodule of $A$.  However, $G=\?M^{-1}K$.
So $K\subset \?MK$.  Hence, by Nakayama's lemma, $K=0$, which is
absurd.

Set $D:=A+G$.  Then $D\supset A$, but $D\neq A$.  Set $D':=\Hom_A(D,A)$.
Then $D'\subset A$.  Also, $\Hom_A(D',A)=D$ because $A$ is Gorenstein.
So $D'\neq A$.  Hence $D'\subset M$.  Therefore $D\supset \Hom_A(M,A)$.
Now, $\Hom_A(M,A)=B$ by Part~(1) of Lemma~\Cs7).
Hence
        $$BI\subset DI=I+GI.$$
 So \Cs8.3) yields $BI\subset I$, which was to be proved.

\dfn9 Let $C/S$ be a flat projective family of integral curves, and
$Y\subset C$ a closed subset.
 Form the subscheme of $\CJ^n$
parameterizing those sheaves $\cI$ such that $\cI^\star$ is invertible
along $Y$, and denote this open subscheme by $\CJ^n_{C,Y}$.

 \prp10
 Let $\vf\:\dr C\to C$ be a birational $S$-map of flat projective
families of integral curves, and $\delta$ its genus change.
  If $\delta=1$, and if $C/S$ is Gorenstein along the conductor
subscheme $Y$, then
  $$\CJ^n_{C,Y}=\CJ^n_{C}-\epsilon_\vf(\CJ^{n-1}_{\dr C}).$$
 \pf We may assume that $S$ is the spectrum of an algebraically closed
field because forming $Y$ commutes with changing $S$ since $C/S$ is
Gorenstein along $Y$ (or since $\delta=1$); see  \Cs1).  Then
$Y$ consists of a single point $Q$.  Apply Lemma~\Cs7); its Part~(3) says
that $\vf$ factors through the map $\psi$ introduced there.  The genus
drop of $\psi$ is 1 by Part~(2) of Lemma~\Cs7).  So $\vf=\psi$ since
$\delta=1$.  Let $\cI$ be a torsion-free rank-1 sheaf on $C$.  Then
$\cI$ is invertible along $Y$ if and only if $\cI^\star$ is so, since
$C$ is Gorenstein along $Y$.  Hence Lemma~\Cs8) yields the
assertion.

\sct4 Filtrations

\dfn1 Let $C/S$ be a flat projective family of integral curves, and
$Y\subset C$ a closed subset.  Fix $n$.  Using the notation of
(2.3) and (3.9), set
        $$\eqalign{W^{n,k}_{Y,\bigomega}
        &:=W^{n,k}_\bigomega\cap(\CJ^{n}_{C,Y}\x\Quot_\bigomega^k),\cr
        \CJ^{n,k}_{C,Y}
        &:=\alpha_\bigomega(W^{n+k,k}_{Y,\bigomega})\subset\CJ^n_C.
        \cr}$$
 Note that $W^{n,k}_{Y,\bigomega}$ and $\CJ^{n,k}_{C,Y}$ are open, the
latter thanks to Theorem (2.4).

Let $\vf\:\dr C\to C$ be a birational $S$-map, and $Y$ now be its
conductor subscheme.  Let $g$ denote the arithmetic genus of $C/S$, set
$m:=2g-2-n$, and use the notation of (3.1) and (3.3).  Finally, form the
following {\it sets}:
  $$\eqalign{
   K^{n,r}_\vf&:=\iota\circ\kappa_\vf(P^{m+r-\delta,r}_\vf)
        \subset\CJ^n_C,\cr
   L^{n,r}_\vf
        &:=(\iota\circ\kappa_\vf)^{-1}
        (\CJ_{C,Y}^{m-\delta+r,\delta-r})\subset P^{n,r}_\vf.\cr}$$
By convention, take $\CJ^{n,-1}_{C,Y}$ and $K^{n,-1}_\vf$ to be empty.
Note that  $L^{n,r}_\vf$ is open.

 \lem2 Let $C$ be a projective integral curve over an algebraically
 closed field.  Let $\vf\:C^\dagger\to C$ be a birational map, $\delta$
its genus drop, and $Y\subset C$ and $\dr Y\subset \dr C$ its conductor
subschemes.  Let $\cN$ be a torsion-free rank-$1$ sheaf on $C$, and set
        $$d:=\length(\cN\cO_{\dr C}/\cN),$$
 where  $\cN\cO_{\dr C}$ stands for the quotient of
$\cN\ox\cO_{\dr C}$ by its torsion subsheaf.
 Then the following statements hold.
 \part1 The natural map $h\:\cN\to \vf_*(\cN\cO_{\dr C})$ is a
presentation of colength $d$.
 \part2 Let $\cL\subset\cN$ be any subsheaf that is invertible along $Y$.
Then
        $$\length(\cN/\cL)\ge\delta-d,$$
 and equality holds if and only if $\cL\cO_{\dr C}=\cN\cO_{\dr C}$.
 \part3 If $\cN\cO_{\dr C}$ is invertible along  $\dr Y$, then
there exists a subsheaf $\cM$ of $\cN$ such that $\cM$ is invertible
along $Y$ and $\cM\cO_{\dr C}=\cN\cO_{\dr C}$.
 \pf
 To prove (1), note that $h$ is injective because $\cN$ is torsion-free
and $\vf$ is birational.  In addition, $\cCok(h)$ is an $\cO_Y$-module
since $\cI_Y$ is the conductor ideal, and so
  $$\cI_Y\vf_*(\cN\cO_{\dr C})=\cI_Y\cN\subset\cN.$$

 To prove (2), note that we have the inclusions,
  $$\cN\cO_{\dr C}\supset\cN\supset \cL
  \hbox{ and }\cN\cO_{\dr C}\supset\cL\cO_{\dr C}\supset \cL.$$
 So since length is additive, we get the equations,
        $$\eqalign{\length(\cN\cO_{\dr C}/\cL)
  &=\length(\cN\cO_{\dr C}/\cN)+\length(\cN/\cL),\cr
   \length(\cN\cO_{\dr C}/\cL)
  &=\length(\cN\cO_{\dr C}/\cL\cO_{\dr C})
    +\length(\cL\cO_{\dr C}/\cL).\cr}$$
 Since $\cL$ is invertible along $Y$, the last term of the second
equation is equal to $\delta$.  Hence
 $$\length(\cN/\cL)=\length(\cN\cO_{\dr C}/\cL\cO_{\dr C})+\delta-d.$$
 The assertion follows immediately.

 To prove (3), first
 take
 $Q\in Y$.  Since $\cN\cO_{\dr C}$ is invertible
along $\vf^{-1}Q$, there is an $a_Q\in\cN_Q$ such that
  $$a_Q(\vf_*\cO_{\dr C})_Q=\cN_Q(\vf_*\cO_{\dr C})_Q.$$
 Now, define $\cM$ to be the subsheaf of $\cN$ such that $\cM_Q$ is
equal to $a_Q\cO_{C,Q}$ if $Q\in Y$, and to $\cN_Q$ if not.  Then $\cM$
has the desired properties. Thus the lemma is proved.

\prp3 Let $C$ be a projective integral curve over an algebraically
closed field.  Let $\vf\:\dr C\to C$ be a birational map, $\delta$ its
 genus drop, $Y\subset C$ and $\dr Y\subset \dr C$ its conductor
subschemes, and $\cI_Y$ its conductor ideal on $C$.  Let $\cK$ be a
torsion-free rank-$1$ sheaf of degree $n$ on $C$, and set
  $$c:=\length(\cK/\cI_Y\cK)\and
    d:=\length(\cK^\star\cO_{\dr C}/\cK^\star),$$
 where  $\cK^\star\cO_{\dr C}$ stands for the quotient of
$\cK^\star\ox\cO_{\dr C}$ by its torsion subsheaf.
 Let $t\in\CJ^n_C$ represent $\cK$.  Then the following statements hold.
 \part1 If $t\in\CJ^{n,k}_{C,Y}$, then $k\ge\delta-d$.  The converse
 holds if $\cK^\star\cO_{\dr C}$ is invertible along $\dr Y$.
 \part2 We have $d\le c$.  In addition, $t\in K^{n,r}_\vf$
if and only if $d\leq r\leq c$.
 \pf We are going to apply Lemma~\Cs2) a number of times with
$\cN:=\cK^\star$.

Let's prove (1).  Suppose $t\in\CJ^{n,k}_{C,Y}$.
 Then $\cK^\star=\cI^\star\otimes\cJ^\star$ where $\cJ\subset\bigomega_C$
 is a subsheaf of colength $k$ and $\cI$ is a torsion-free rank-1 sheaf
on $C$ such that $\cI^\star$ is invertible along $Y$ and where
$\cJ^\star$ is not.  The inclusion $\cJ\into\bigomega_C$ yields an
injection $\cO_C\into\cJ^\star$.  Tensoring it with $\cI^\star$, we
obtain an injection $\cI^\star\into\cK^\star$.
  Form the subsheaf $\cL\subset\cK^\star$ that is equal to
$\cI^\star$
 along
$Y$ and to $\cK^\star$ off $Y$.
Then we have these relations,
  $$k=\length(\cJ^\star/\cO_C)\geq\length(\cK^\star/\cL)\geq
  \delta-d,$$
 since the operator $^\star$ preserves colength, since $\cI^\star$ is
invertible along $Y$, and since Part (2) of Lemma~\Cs2) applies.

Conversely, suppose $\cK^\star\cO_{\dr C}$ is invertible along $\dr Y$.
By Part (3) of Lemma~\Cs2), there is a subsheaf
$\cM\subset\cK^\star$ such that $\cM$ is invertible along $Y$ and
$\cM\cO_{\dr C}=\cK^\star\cO_{\dr C}$.  By Part (2) of Lemma~\Cs2),
        $$\length(\cK^\star/\cM)= \delta-d.$$

In addition, suppose $k\geq\delta-d$.  Set
$r:=k-\delta+d$, take a simple closed point $R\in C$, and
set $\cI:=\cM^\star(rR)$.  Then $\cI^\star$ is invertible along
$Y\cup\{R\}$, and
        $$\eqalign{\deg\cI
        &=\deg\cK+\length(\cK^\star/\cM)+\length(\cM/\cI^\star)\cr
        &=n+\delta-d+r=n+k.\cr}\eqno\Cs3.1)$$

Set $\cH:=\cHom(\cI^\star,\cK^\star)$ and embed $\cO_C$ in $\cH$ as the
subsheaf of multiples of the inclusion map $\cI^\star\into\cK^\star$.
Form the subsheaf $\cL$ of $\cH$ whose stalk $\cL_Q$ is $\cH_Q$ if $Q\in
Y\cup\{R\}$, and is $\cO_{C,Q}$ otherwise.  Because of the latter
condition, $\cL$ is invertible where $\cI^\star$ is not.  Set
$\cJ:=\cL^\star$.  Consider the natural map from $\cI^\star\ox\cJ^\star$ into
$\cK^\star$; it is an isomorphism, because $\cI^\star$ is invertible along
$Y\cup\{R\}$, and because off it, $\cI^\star$ is equal to $\cK^\star$ and
$\cJ^\star$ to $\cO_C$.

There is a natural injection of $\cO_C$ into $\cJ^\star$.  Its colength
is equal to the colength of $\cI^\star$ in $\cK^\star$.  The latter is
equal to $\deg\cI-\deg\cK$, so to $k$ by \Cs3.1).  Applying $^\star$
yields an injection of $\cJ$ into $\bigomega_C$, and its cokernel is
also of length $k$.  In sum, the pair $(\cI,\cJ)$ is represented by a
point of $W^{n+k,k}_{\bigomega,Y}$, and its image under
$\alpha_\bigomega$ is $t$.  Thus (1) holds.

Let's prove (2).  Consider the natural map $h\:\cK^\star\to
\vf_*(\cK^\star\cO_{\dr C})$; it is a presentation of colength $d$ by
Part~(1) of Lemma~\Cs2).  Set $\dr\cL:=(\cK^\star\cO_{\dr C})^\star$,
and let $h^\star\:\vf_*\dr\cL\into\cK$ be the copresentation associated
to $h$.  Since $\cI_Y\cK$ is an $\vf_*\cO_{\dr C}$-module contained in
$h^\star(\vf_*\dr\cL)$, there is an $\cO_{\dr C}$-submodule
$\dr\cI\subset\dr\cL$ such that $h^\star(\vf_*\dr\cI)=\cI_Y\cK$.  The
quotient $\dr\cL/\dr\cI$ is an $\cO_{\dr Y}$-module.  Furthermore,
        $$\eqalign{\length(\dr\cL/\dr\cI)
  &=\length(h^\star(\vf_*\dr\cL)/h^\star(\vf_*\dr\cI))\cr
  &=\length(\cK/h^\star(\vf_*\dr\cI))
        -\length(\cK/h^\star(\vf_*\dr\cL))\cr
  &=\length(\cK/\cI_Y\cK)
        -\length((\vf_*\dr\cL)^\star/\cK^\star)=c-d.\cr}$$
 Thus $0\le c-d$.

Suppose $t\in K^{n,r}_\vf$.  Then there is a presentation $u\:\cK^\star\to
\vf_*\dr\cJ$ of colength $r$.  Now, $u$ factors through $h\:\cK^\star\to
\vf_*(\cK^\star\cO_{\dr C})$.  So $d\le r$.  Consider the associated
copresentation $u^\star\: \vf_*(\dr\cJ)^\star\to\cK$.  Since
$\cCok(u^\star)$ is an $\cO_{Y}$-module of length $r$, we have $r\leq
c$, as desired.

Conversely, suppose $d\leq r\leq c$.  Then there is  a
subsheaf $\dr\cK$ of $\dr\cL$ with $\dr\cK\supset\dr\cI$ and
$\length(\dr\cL/\dr\cK)=r-d$.  So the composition,
  $$\vf_*\dr\cK\to\vf_*\dr\cL@>h^\star>>\cK$$
 is a copresentation of colength $r$ precisely.  So $t\in K^{n,r}_\vf$,
and the proposition is proved.

\prp4 Let $C/S$ be a flat projective family of integral curves, and
$Y\subseteq C$ a closed subset.  Fix $n$.  Then the following statements
hold.
 \part 1 The open subsets $\CJ^{n,k}_{C,Y}$ of $\CJ^n_C$
provide an exhaustive increasing filtration:
        $$\CJ^n_{C,Y}=\CJ^{n,0}_{C,Y}\subseteq\cdots\textstyle
        \subseteq\CJ^{n,k}_{C,Y}\subseteq\cdots\subseteq
        \textstyle\bigcup_{k\ge 0}\CJ^{n,k}_{C,Y}=\CJ^n_C.$$
 \part 2 If, as $s$ ranges over the geometric points of $S$, there is an
upper bound $\Delta$ on the genus drop to the partial normalization
along $Y(s)$ of $C(s)$, then
        $$\CJ^{n,\Delta}_{C,Y}=\CJ^n_C.$$
  \pf It follows directly from the definitions
 that $\alpha^{n,0}_\bigomega$ is equal to the identity of $\CJ^n_C$;
whence, $\CJ^n_{C,Y}=\CJ^{n,0}_{C,Y}$.  Now, forming $\CJ^{n+k}_{C,Y}$
commutes with changing $S$; whence, so does forming
$\CJ^{n,k}_{C,Y}$.  Hence, it suffices
 to check on the geometric fibers that
$\CJ^{n,k}_{C,Y}\subseteq\CJ^{n,k+1}_{C,Y}$ for any $k$ and that
$\CJ^{n,\Delta}_{C,Y}=\CJ^n_C$.  Thus we may
 assume that $S$ is the spectrum of an algebraically closed field.

Let $t\in\CJ^n_C$ be a closed point, and $\cK$ a corresponding sheaf.
Let $\vf\:\dr C\to C$ be the partial normalization of $C$ along $Y$, and
set $d:=\length(\cK^\star\cO_{\dr C}/\cK^\star)$.
Now, $\cK^\star\cO_{\dr C}$
is invertible along the conductor subscheme because $\dr C$ is smooth
there.  So Part (1) of Proposition~\Cs3) applies, and it says that
$t\in\CJ^{n,k}_{C,Y}$ if and only if $k\ge\delta-d$.  Hence, if
$t\in\CJ^{n,k}_{C,Y}$, then $t\in\CJ^{n,k+1}_{C,Y}$.  In addition,
$t\in\CJ^{n,\Delta}_{C,Y}$ because $\Delta\geq\delta$.  Thus the
proposition is proved.

 \prp5 Let $C/S$ be a flat projective family of integral curves,
 $\vf\:\dr C\to C$ a birational $S$-map, $\delta$ its genus drop,
$Y\subset C$ and $\dr Y\subset\dr C$
its conductor subschemes.  Fix $n$.
 \part 1 Then $\bigcup_{r\ge 0}K_\vf^{n,r}=\CJ^n_C$.  If the genus drop
 to
 the partial normalization of every geometric
fiber $C(s)$ along $Y(s)$ is bounded
by $\Delta$, then $\bigcup_{r\ge 0}^\Delta K_\vf^{n,r}=\CJ^n_C$
 \part 2 If $\dr C/S$ is smooth along $\dr Y$, then
$\bigcup_{r\ge0}^\delta K_\vf^{n,r}=\CJ^n_C$.
 \part 3 If $\dr C/S$ is smooth along $\dr Y$ and if $C/S$ is
Gorenstein along $Y$, then $K_\vf^{n,\delta}=\CJ^n_C$.
 \pf Let $t$ be a geometric point of $\CJ^n_C$, and $\cK$ a
corresponding sheaf on $C(t)$.  Set $\dr\cK:=\cK\cO_{\dr C(t)}$ and
$\dr\cL:=\cK^\star\cO_{\dr C(t)}$, both of these products being the
quotient of the corresponding tensor product by its torsion subsheaf.
Let $\psi\:\wt C\to C(t)$ be the partial normalization of $C(t)$ along
$Y(t)$, and $\wt\delta$ the genus drop along $\psi$.  Similarly, set
$\wt\cL:=\cK^\star\cO_{\wt C}$.

Apply Lemma~\Cs2) with $\cN:=\cK^\star$.  Since $Y(t)$ contains the
conductor subscheme of $\vf(t)$, Part (1) of Lemma~\Cs2) implies that
the natural map $h\:\cK^\star\to \vf(t)_*\dr\cL$ is a presentation for
$\vf$, say of colength $d$.  Hence $t\in K_\vf^{n,d}$.

Let $u\:\cK^\star\to\psi_*\wt\cL$ be the natural injection, $r$ its
colength.  Then $u$ factors through $h$ because $\psi$ factors through
$\vf(t)$. So $d\le r$. Now, $\wt\cL$ is invertible along $\psi^{-1}Y(t)$
because $\wt C$ is smooth there. So Parts (2) and (3) of Lemma~\Cs2)
imply that $0\le\wt\delta-r$.  Hence $d\le\wt\delta$.  If
$\wt\delta\leq\Delta$, then $d\leq\Delta$.  Thus (1) holds.

Suppose $\dr C/S$ is smooth along $\dr Y$.  Then $\dr C(s)$ is the
partial normalization of $C(s)$ along $Y(s)$ for every geometric point
$s$ of $S$.  Taking $\Delta:=\delta$ in (1), we obtain (2).

Since $\dr C(t)$ is smooth along $\dr Y(t)$, we have $\wt\delta=\delta$.
 So
 $d\le\delta$.  In addition, $\dr\cK$ is also invertible along
$\dr Y(t)$.  Apply Lemma~\Cs2) to $\cN:=\cK$ now.
Parts (2) and (3) imply that
        $$\length(\dr\cK/\cK)\le\delta.\eqno\Cs5.1)$$

Let $\cI$ be the conductor ideal on $C(t)$ of $\vf(t)$.  Set
  $$c:=\length(\cK/\cI\cK) \and \delta':=\length(\cO_{C(t)}/\cI).$$
 Since $\cI=\cI\cO_{\dr C(t)}$ and since $\dr\cK$ is invertible along
$\dr Y(t)$, we have
        $$\length(\dr\cK/\cI\cK)=\delta+\delta'.\eqno\Cs5.2)$$
{}From \Cs5.1) and \Cs5.2), we get
        $$c=\length(\dr\cK/\cI\cK)-
        \length(\dr\cK/\cK)\ge\delta'.$$

Suppose now that $C/S$ is Gorenstein along $Y$. Then $\delta'=\delta$,
and so $d\leq\delta\leq c$.  Hence, Part (2) of Proposition~\Cs3)
implies that $t\in K^{n,\delta}_\vf$. Thus the proposition is proved.

 \rmk6 Let $C$ be a projective integral curve of arithmetic genus $g$
over an algebraically closed field, and $Y$ a closed point.  Let
$\vf\:\dr C\to C$ be the partial normalization of $C$ at $Y$, and
$\delta$ the genus drop.  Fix $n$.  In this setting, the following
statement holds:
 $$\hbox{\it If $K_\vf^{n,\delta}=\CJ^n_C$, then $C$ is Gorenstein at
 $Y$.}$$

To prove it, let $t\in\CJ^n_C$ be a closed point, and $\cK$ a
corresponding sheaf.  Say $t\in K_\vf^{n,r}$.  Then there is a
presentation $h\:\cK^\star\to \vf_*\dr\cL$ of colength $r$.  Let
$h^\star\:\vf_*{\dr\cL}^\star\to\cK$ be the associated copresentation.
Let $\cI_Y$ be the conductor ideal.  Then
$\cI_Y\cK\subset\cIm(h^\star)$.  Hence
        $$\length(\cK/\cI_Y\cK)\ge\colength(h^\star)=r.$$

Assume now that $\cK$ is invertible at $Y$.  Then
        $$\length(\cK/\cI_Y\cK)=\length(\cO_C/\cI_Y).$$
 However,
 $\length(\cO_C/\cI_Y)\le\delta$,
 and equality holds (if and)
only if $C$ is Gorenstein at $Y$.  Hence, $r\le\delta$, and if
$r=\delta$, then $C$ is Gorenstein at $Y$.  Finally, if
$K_\vf^{n,\delta}=\CJ^n_C$, then all the preceding considerations apply
to any closed point $t\in J^n_C$, and so $C$ is Gorenstein at $Y$.

In any event, $K_\vf^{n,\delta}$ is closed and contains $\iota
J^{2g-2-n}_C$, which is open.  Hence, if $\CJ^n_C$ is irreducible, then
$K_\vf^{n,\delta}=\CJ^n_C$.  Thus we obtain the following corollary:
  $$\hbox{\it If $\CJ^n_C$ is irreducible, then $C$ is Gorenstein.}$$
  In fact, if $\CJ^n_C$ is irreducible, then $C$ has embedding dimension
2 at every singular point; this stronger result was obtained in
\cite{KK81, (1), \p.277} and \cite{R80, Theorem A, \p.212}.

On the other hand, our proof yields the following sharper result:
  $$\hbox{\it If the point representing $\bigomega_C$ lies in the
closure of $J^{2g-2}_C$, then $C$ is Gorenstein.}$$
 Indeed, if so, then $0\in J^0_C$ lies in the closure of $\iota
J^{2g-2}_C$.  So $0\in K_\vf^{0,\delta}$.  Hence, if we take $t:=0$,
$n:=0$ and $r:=\delta$, then the considerations above imply that $C$ is
Gorenstein.

\eg7 Here is an example where $\kappa^{n,r}_\vf$ is not surjective for
any integer $r\geq 0$.  Let $C$ be a projective integral curve over an
algebraically closed field.  Assume that $C$ has a unique singularity
$Q$.  Assume also that the normalization map $\vf\:\dr C\to C$ has genus
drop $\delta$ greater than 1, and that the conductor subscheme $Y\subset
C$ is of length 1, hence equal to $Q$.  A concrete choice is
     $$C:=\{(s^5:s^2t^3:st^4:t^5)\in\IP^3 \mid (s:t)\in\IP^1\}\and
          Q:=(1:0:0:0).$$

 Since $\dr C$ is smooth, $\CJ^n_{\dr C}$ is irreducible. Since $Y=Q$,
Proposition~(3.5) implies that $P_\vf^{n,r}$ is a Grassmann bundle over
$\CJ^n_{\dr C}$.  Hence, $P_\vf^{n,r}$ is irreducible as well for every
$r\geq 0$.  If $\kappa^{n,r}_\vf$ were surjective for some $r$, then
$\CJ^{n+\delta-r}_C$ would be irreducible.  However, $C$ is not
Gorenstein at $Q$ because $\delta>1$ and $\length(\cO_Y)=1$.  So
$\CJ^{n+\delta-r}_C$ is reducible by \cite{KK81, (1), \p.277};
alternatively, see Remark~\Cs6).  Hence $\kappa^{n,r}_\vf$ is not
surjective for any $r\geq 0$.

 \prp8 Let $C/S$ be a flat projective family of integral curves,
 $\vf\:\dr C\to C$ a birational $S$-map, and $\dr Y\subset\dr C$ the
conductor subscheme.  If $\dr C/S$ is smooth along $\dr Y$, then
$L^{n,r}_\vf$ is equal to the open subscheme of $P^{n,r}_\vf$
parameterizing presentations whose adjoint is surjective.
 \pf Forming the adjoint commutes with changing the base; moreover,
surjectivity may be checked on the fibers.  Hence it suffices to prove
the following statement: let $t$ be a geometric point of $P_\vf^{n,r}$,
and $u\:\cI\to \vf(t)_*\dr\cJ$ a corresponding presentation on the fiber
$C(t)$; then $t$ lies in $L^{n,r}_\vf$ if and only if the adjoint
$u^\sharp\:\vf(t)^*\cI\to\dr\cJ$
 is surjective.

  Form $\dr\cI:=\cI\cO_{\dr C(t)}$, the quotient of $\cI\ox\cO_{\dr
C(t)}$ by its torsion subsheaf.  Form the natural presentation
$h\:\cI\to \vf(t)_*\dr\cI$, and let $d$ be its colength.
  Then $u$ factors as follows:
        $$u\:\cI@>h>>\vf(t)_*\dr\cI@>\vf(t)_*\dr u>> \vf(t)_*\dr\cJ$$
 where $\dr u\:\dr\cI\to\dr\cJ$ is the induced map.  Then $\vf(t)_*\dr
u$ is injective.  Now, $u$ is of colength $r$.  Hence $\vf(t)_*\dr u$ is
an isomorphism if and only if $d=r$, or equivalently, $d\ge r$.  Now,
$u^\sharp$ is surjective if and only if $\dr u$ is an isomorphism.
Therefore, $u^\sharp$ is surjective if and only if  $d\ge r$.

On the other hand, $t$ lies in $L^{n,r}_\vf$ if and only if the point
representing $\cI^\star$ lies in $\CJ^{m,\delta-r}_{C,Y}$, where
$\delta$ is the genus drop along $\vf$ and $m:=2g-2+r-\delta-n$, where
$g$ is the arithmetic genus of $C/S$.
  Let's apply Part~(1) of Proposition~\Cs3) with
$\cK:=\cI^\star$.  We conclude that $t$ lies in $L^{n,r}_\vf$ if and
only if $\delta-r\geq\delta-d$.  Therefore, by the preceding paragraph,
$t$ lies in $L^{n,r}_\vf$ if and only if $u^\sharp$ is surjective.  Thus
the proposition is proved.

 \thm9 Let $C/S$ be a flat projective family of integral curves,
$\vf\:\dr C\to C$ be a birational $S$-map, $\delta$ the genus drop, and
$Y\subset C$ the conductor subscheme.  Fix $n$.  Let $g$ be the
arithmetic genus of $C/S$, and set $m:=2g-2-n$.  Fix $k$ such that $0\le
k\le\delta$.  Then the following statements hold.
 \part 1 If $\dr C/S$ is smooth along the conductor subscheme $\dr Y$,
then $\iota\circ\kappa_\vf$ restricts to a
 closed embedding of schemes, $L_\vf^{m-k,\delta-k}\into\CJ^{n,k}_{C,Y}$.
 \part 2 If $\dr C/S$ is smooth along $\dr Y$ and $C/S$ is
 Gorenstein along $Y$, then set-theoretically
        $$\eqalign{\CJ^n_C-\CJ^{n,k-1}_{C,Y}&=K_\vf^{n,\delta-k},\cr
        \CJ^{n,k}_{C,Y}-\CJ^{n,k-1}_{C,Y}&=\iota\circ
        \kappa_\vf(L_\vf^{m-k,\delta-k}).\cr}$$
 \pf To prove (1), we may change the base via an \'etale
surjective map $S'\to S$.  We may thus assume that there exists a
universal presentation on $C\x P^{m-k,\delta-k}_\vf$.
Since $\kappa_\vf$ is proper and $\iota$ is an
isomorphism, it is enough to show that the restriction,
        $$\kappa_\vf|L^{m-k,\delta-k}_\vf,$$
 is a monomorphism.

 Let $h\:\cI\to \vf_{T*}\dr\cJ$ be a presentation corresponding to a
$T$-point of $L^{m-k,\delta-k}_\vf$.  By Proposition~\Cs8), the adjoint
$h^\sharp\: \vf_T^*\cI\to\dr\cJ$ is surjective.  Set $U:=C-Y$.  Then
$\vf_T$ is an isomorphism over $U_T$.  Let $i\: U_T\to C_T$ denote the
inclusion, and $j\: U_T\to\dr C_T$ its lifting.
Let $\dr\cI$ denote the image of the natural map $u\: \vf_T^*\cI\to
j_*i^*\cI$, and $v\:\cI\to \vf_{T*}\dr\cI$ its adjoint.
  We need only show that $v$ is a presentation equivalent to $h$.

Consider the following commutative diagram on $\dr C_T$:
 $$\CD
 \vf_T^*\cI  @>h^\sharp>>    \dr\cJ\\
   @VVuV                @VVwV\\
 j_*i^*\cI @>j_*j^*h^\sharp>> j_*j^*\dr\cJ
 \endCD$$
 where $w$ is the natural map.  Since $\cCok(h)$ is supported
 outside $U_T$, the map $j_*j^*h^\sharp$ is an isomorphism.  Since the
associated points of $\dr\cJ$ lie in $j(U_T)$, the map $w$ is injective.
Since $h^\sharp$ is surjective, and $\dr\cI$ is the image of $u$, the
map $j_*j^*h^\sharp$ induces an isomorphism $c\:\dr\cI\risom\dr\cJ$
such that $h^\sharp=cu$.  By adjunction, $h$ is equivalent to $v$.  Thus
(1) is proved.

Consider (2) now.
 The asserted equations may be checked on the fibers over $S$.  Since
$C/S$ is Gorenstein along $Y$, forming the presentation scheme commutes
with changing $S$; see (3.1).
 So
 we may assume that $S$ is the spectrum of an algebraically closed
field.  Let $t$ be a closed point of $\CJ^n_C$, and $\cK$ a
corresponding sheaf on $C$.  Set
        $$c:=\length(\cK/\cI_Y\cK)\and d:=
        \length(\cK^\star\cO_{\dr C}/\cK^\star).$$

Since $\dr C$ is smooth along
$\dr Y$, Part~(1) of Proposition~\Cs3) says that
        $$t\in\CJ^{n,k-1}_{C,Y}{\rm \  if\  and\  only\  if\  }
        k-1\geq\delta-d.\eqno\Cs9.1)$$
 In addition, since $C/S$ is Gorenstein along $Y$,
Part~(3) of Proposition~\Cs5) says that
there is a presentation $u\:\cK^\star\to \vf_*\dr\cJ$ of colength $\delta$.
Since the associated copresentation $u^\star\: \vf_*{\dr\cJ}^\star\to\cK$ is
of colength $\delta$ as well, we have $c\geq\delta$. Hence, Part (2) of
Proposition~\Cs3) says that
        $$t\in K^{n,\delta-k}_\vf{\rm \  if\  and\  only\  if\  }
        d\leq\delta-k.\eqno\Cs9.2)$$
 Since either $d\geq\delta-k+1$ or $d\leq\delta-k$,
it follows from \Cs9.1) and \Cs9.2) that
$\CJ^{n,k-1}_{C,Y}$ and $K^{n,\delta-k}_\vf$ are disjoint and
cover $\CJ^n_C$.
 Thus the first asserted equation holds.

By Part~(1) of Proposition~\Cs3), we have
        $$t\in\CJ^{n,k}_{C,Y}-\CJ^{n,k-1}_{C,Y}
        {\rm \  if\  and\  only\  if\  }k=\delta-d.\eqno\Cs9.3)$$
 On the other hand, by definition,
        $$\iota\circ\kappa_\vf(L^{m-k,\delta-k}_\vf)=
        \CJ^{n,k}_{C,Y}\cap K^{n,\delta-k}_\vf.\eqno\Cs9.4)$$
 By Part (2) of Proposition~\Cs3), $t\in K^{n,\delta-k}_\vf$ if and only if
$d\leq\delta-k$. By Part (1) of Proposition~\Cs3),
$t\in\CJ^{n,k}_{C,Y}$ if and only if $k\geq\delta-d$.  Hence, \Cs9.4)
implies that
        $$t\in\iota\circ\kappa_\vf(L^{m-k,\delta-k}_\vf)
        {\rm \  if\  and\  only\  if\  }k=\delta-d.\eqno\Cs9.5)$$
 Combining \Cs9.3) and \Cs9.5), we obtain the second asserted equation.
Thus the Theorem is proved.

\sct5 The lifted Abel map

 \prp1 Let $\cF$ be a torsion-free rank-$1$ sheaf on a flat projective
family $C/S$ of integral curves.  Let $\vf\:\dr C\to C$ be a birational
 $S$-map,  $Y\subseteq C$  its conductor subscheme, and $\cI_Y$ its
conductor ideal.
  Then the following statements hold.
 \part1 There is, up to equivalence, at most one presentation
$h\:\cF\to\vf_*\dr\cF$ whose adjoint $h^\sharp\:\vf^*\cF\to\dr\cF$ is
surjective.
 \part2 Assume that there is a presentation $h$ as in {\rm (1)}, and
let $h^\star\:\vf_*{\dr\cF}^\star\to\cF^\star$ be the associated map.
If $\cIm(h^\star)=\cI_Y\cF^\star$, then $h$ is the only presentation
whose source is $\cF$.
 \part3 Assume that  $Y$ is flat, and that forming $Y$
commutes with changing $S$.  Assume that $\cF^\star$ is invertible along
$Y$.  Then there is one and only one presentation $h\:\cF\to\vf_*\dr\cF$
whose source is $\cF$; its
  adjoint $h^\sharp\:\vf^*\cF\to\dr\cF$ is surjective, and
$\cIm(h^\star)=\cI_Y\cF^\star$.
 \pf Let's prove (1) and (2) first.  Let $u\:\cF\to\vf_*\dr\cG$ be any
presentation, and $h$ be as in (1).  Since $h$ is a presentation,
$\cKer(h^\sharp)$ contains no fiber of $C/S$ in its support.  So
$u^\sharp(\cKer(h^\sharp))$ is a subsheaf of $\dr\cG$ whose support
contains no generic point of any fiber.  Hence, since $\dr\cG$ is
torsion-free, $u^\sharp(\cKer(h^\sharp))=0$.
  So $\cKer(h^\sharp)\subset
\cKer(u^\sharp)$.  Hence, since $h^\sharp$ is surjective, $u^\sharp$
factors through a unique map $v\:\dr\cF\to\dr\cG$.  By adjunction,
$u=\vf_*v\circ h$.

To prove (1), assume that $u^\sharp$ is surjective as well.  By
reversing the roles of $u$ and $h$ in the above argument, we get an
inverse to $v$.  So $u$ and $h$ are equivalent.
 Thus (1) holds.

To prove (2), assume that $\cIm(h^\star)=\cI_Y\cF^\star$.  Let
$u^\star\:\vf_*{\dr\cG}^\star\to\cF^\star$ and $v^\star\:{\dr\cG}^\star
\to{\dr\cF}^\star$ be the maps associated to $u$ and $v$.
 They are injective.  Moreover, $u^\star=h^\star\circ\vf_*v^\star$;
whence, $\cIm(h^\star)\supset\cIm(u^\star)$.  Since $u^\star$ is a
copresentation, $\cIm(u^\star)\supset\cI_Y\cF^\star$.  So
$\cIm(u^\star)=\cIm(h^\star)$.  Therefore, $\vf_*v^\star$ is an
isomorphism.  By Lemma~(3.4), $v^\star$ is then an isomorphism, and so
$v$ is an isomorphism.  Hence $u$ and $h$ are equivalent.  Thus (2)
holds.

Let's prove (3).
  It follows from the short exact sequence,
        $$0\to\cO_C/\cI_Y\to \vf_*\cO_{\dr C}/\cI_Y\to
        \vf_*\cO_{\dr C}/\cO_C\to0,$$
 that the conductor subscheme $\dr Y\subset\dr C$ is flat, and that
forming $\dr Y$ commutes with changing $S$.  Since $\cF^\star$ is
invertible along $Y$, the pullback $\vf^*(\cF^\star)$ is
 invertible along $\dr Y$.  Hence the restrictions $\cF^\star|Y$ and
$\vf^*(\cF^\star)|\dr Y$ are flat, and forming them commutes with
changing $S$.  Therefore, $\cI_Y\cF^\star$ and $\cI_{\dr Y}\vf^*(
\cF^\star)$ are flat, and forming them commutes with changing $S$.

Since $\cI_Y$ and $\cI_{\dr Y}$ are the conductor ideals,
$\cI_Y\cF^\star = \vf_*(\cI_{\dr Y}\vf^*(\cF^\star))$.  So the inclusion
$\cI_Y\cF^\star\to\cF^\star$ is a copresentation.  Consider the
associated presentation,
        $$h\:\cF\to \vf_*\dr\cF\  \where\
        \dr\cF:=(\cI_{\dr Y}\vf^*(\cF^\star))^\star.$$
 By construction, $\cIm(h^\star)=\cI_Y\cF^\star$.  Next, we'll prove
that the adjoint $h^\sharp\:\vf^*\cF\to\dr\cF$ is surjective.  Then (2)
will imply that $h$ is unique.  Thus $h$ will be as asserted.

Let $s\in S$ and set $\dr\cJ:=\cIm(h^\sharp(s))$.  Then
$h(s)$ factors through $\vf(s)_*\dr\cJ$, and the associated
copresentation factors as follows:
        $$h(s)^\star\: \vf(s)_*\dr\cF(s)^\star
        @>v>> \vf(s)_*{\dr\cJ}^\star @>w>> \cF(s)^\star.$$
 So $\cIm(h(s)^\star)\subset\cIm(w)$.  Now, $\cIm(h(s)^\star)$ is equal
to $\cI_{Y(s)}\cF(s)^\star$, and $\cI_{Y(s)}\cF(s)^\star$ is the largest
$\cO_{\dr C(s)}$-module contained in $\cF(s)^\star$ because
$\cF(s)^\star$ is invertible along $Y(s)$.  Hence
$\cIm(h(s)^\star)=\cIm(w)$.  Hence $v$ is bijective.  So $v^\star\:
\vf(s)_*\dr\cJ\to\vf(s)_*\dr\cF(s)$ is bijective.  Hence
$\dr\cJ=\dr\cF(s)$.  Hence $h^\sharp(s)$ is surjective for each $s\in
S$.  Therefore, $h^\sharp$ is surjective.  The proof is now complete.

 \eg2 Preserve the conditions of \Cs1).  Suppose $Y$ is flat and
forming $Y$ commutes with changing $S$.  (As noted in (3.1), these
conditions hold if $\delta=1$ or if $C/S$ is Gorenstein along $Y$.)
Since $\bigomega_C^\star=\cO_C$, by Part (3) of Proposition~\Cs1), there
is a unique presentation $h\:\bigomega_C\to\vf_*\dr\cF$.  In this case,
 the proof shows that $\dr\cF=\cI_{\dr Y}^\star$.  If also
 $C/S$ is Gorenstein along $Y$, then the
comorphism $\vf^\sharp\:\cO_C\to\vf_*\cO_{\dr C}$ is the unique
presentation with source $\cO_C$, again by Part (3) of
Proposition~\Cs1).

 Though important, the sheaves $\bigomega_C$ and $\cO_C$ are not the
only examples of a sheaf $\cF$ admitting one and only one presentation
with $\cF$ as source.  For instance, over an algebraically closed field,
consider the map,
        $$\vf\:\IP^1\to\IP^4\hbox{ given by }
        \vf(s:t):=(s^7:s^3t^4:s^2t^5:st^6:t^7).$$
 Set $C:=\vf(\IP^1)$ and $Q:=\vf(1:0)$.  Define an ideal
$\cF\subset\cO_C$ by $\cF_Q:=(t^4,t^5)\cO_{C,Q}$ and $\cF_R:=\cO_{C,R}$
for $R\in C-Q$.  Then $\cF\cO_{\IP^1}=t^4\cO_{\IP^1}$, which is the
conductor ideal.  If there were a presentation of the form
$\cF\to\vf_*(t^n\cO_{\IP^1})$ for $n\le 3$, then we'd have
$t^7\in\cF_Q$, a contradiction.  So $\cF\to\vf_*(t^4\cO_{\IP^1})$ is the
unique presentation with $\cF$ as source.  Note that $\cF_Q$ is not
principal.  Neither is $\cF^\star_Q$ because
$\bigomega_{C,Q}\cong(t^4,t^5,t^6)\cO_{C,Q}$.

 \dfn3 Let $C/S$ be a flat projective family of integral curves of
arithmetic genus $g$.  Let $\vf\:\dr C\to C$ be a birational $S$-map,
$\delta$ its genus change, $Y\subset C$ and $\dr Y\subset \dr C$ its
conductor subschemes; see
(3.1).  As in (3.9), let $\CJ^m_{C,Y}$ denote
the open subscheme of $\CJ^m_C$ of sheaves $\cI$ such that $\cI^\star$
is invertible along $Y$.  Let $\cF$ and $\dr\cF$ be torsion-free rank-$1$
sheaves on $C/S$ and $\dr C/S$ of degree $f$ and $f'$, respectively.  As
in (4.1), set
        $$\eqalign{W^{m,k}_{Y,\cF}
        :=&W^{m,k}_\cF\cap(\CJ^m_{C,Y}\x\Quot^k_\cF),\cr
        W^{m,k}_{\dr Y,\dr\cF}
        :=&W^{m,k}_{\dr\cF}\cap(\CJ^m_{\dr C,\dr Y}\x\Quot^k_{\dr\cF}).\cr}$$
 For convenience, set $W^{m,k}_{Y,C}:=W^{m,k}_{Y,\cO_C}$ and
$W^{m,k}_{\dr Y,\dr C}:=W^{m,k}_{\dr Y,\cO_{\dr C}}$.

Assume that there is a presentation,
        $$h\:\cF\to\vf_{*}\dr\cF,$$
 whose adjoint $h^\sharp\:\vf^*\cF\to\dr\cF$ is surjective.  Define an
$S$-map,
        $$\theta_h
\:\Quot^1_{\dr\cF}\to\Quot^1_\cF,$$
 as follows:  given a $T$-point $q$ of $\Quot^1_{\dr\cF}$, let
$\dr\cJ\subset\dr\cF_T$ be the corresponding subsheaf, and define
$\theta_h(q)$ to be the $T$-point of $\Quot^1_\cF$ representing the
subsheaf $\cJ\subset\cF_T$ given by
        $$\cJ:=h_T^{-1}(\vf_{T*}\dr\cJ)\subset\cF_T.$$

To see that $\theta_h$ is well defined, we have to see that $\cF_T/\cJ$
is flat and of length 1.  To do so, it is enough to see that the
following composite map $u$ is surjective:
        $$u\:\cF_T@>h_T>>\vf_{T*}\dr\cF_T\to\vf_{T*}(\dr\cF_T/\dr\cJ).$$
 So suppose there is a point $t\in T$ such that the fiber $u(t)$ is
not surjective.  Since $\dr\cF(t)/\dr\cJ(t)$ is of length 1, it follows
that $h(t)$ factors through $\vf(t)_*\dr\cJ(t)$, and thus
the adjoint $h^\sharp(t)$ factors through $\dr\cJ(t)$.  However,
 on the contrary,
$h^\sharp$ is surjective.  Thus $u$ is surjective,
and so $\theta_h$ is well defined.

Next we are going to define an $S$-map,
        $$\rho^m\:\  \CJ^m_{C,Y}\to\CJ^{m-2\delta}_{\dr C}.$$
 Let $\cI$ be a torsion-free rank-1 sheaf on
$C_T/T$ such that $\cI^\star$ is invertible along $Y_T$.  Then
the pullback $\vf_T^*(\cI^\star)$ is torsion-free rank-1 on $\dr C_T/T$.
Sending  $\cI$ to $\vf_T^*(\cI^\star)^\star$ defines $\rho^m$.  We will
often abbreviate $\rho^m$ by $\rho$.

Finally, form the following open subscheme:
        $$V^m_{\vf,\dr\cF}
        :=(\rho\x1)^{-1}W^{m-2\delta,1}_{\dr Y,\dr\cF}
        \subset\CJ^m_{C,Y}\x\Quot^1_{\dr\cF}.$$
 It represents the \'etale sheaf associated to the subfunctor whose
$T$-points are those pairs $(\cI,\dr\cJ)$ such that $\cI^\star$ is
invertible along $Y_T$ and $\vf_T^*(\cI^\star)$ is invertible where
${\dr\cJ}^\star$ is not invertible.
  By sending $(\cI,\dr\cJ)$ to
        $$\bigl(\cI^\star\ox(h_T|\cJ)^\star\bigr)^\star,\where
        \cJ:=h_T^{-1}(\vf_{T*}\dr\cJ)\and
        (h_T|\cJ)\:\cJ\to\vf_{T*}\dr\cJ,$$
 we define an $S$-map,
        $$\lambda^m_{\vf,h}\:V^m_{\vf,\dr\cF}\to
        P^{r',\delta'}_\vf \where r':=m-1+f'-2g+2
        \and \delta':=\delta+f'-f.$$
 We will call $\lambda^m_{\vf,h}$ the {\it lifted Abel map\/} of degree
$m$.   We will often abbreviate
$\lambda^m_{\vf,h}$ by $\lambda_{\vf,h}$.

To see that $\lambda_{\vf,h}$ is well defined, we have to see that
$\cI^\star\ox(h_T|\cJ)^\star$ is a copresentation of colength
$\delta'$.  Consider the following commutative diagram,
        $$\CD
        \cJ          @>h_T|\cJ>>    \vf_{T*}\dr\cJ    @>>> \cCok(h_T|\cJ)\\
         @VVV                           @VVV                @VvVV \\
        \cF_T        @>h_T>>       \vf_{T*}\dr\cF_T  @>>>  \cCok(h_T)\\
           @VuVV                        @VVV \\
        \vf_{T*}(\dr\cF_T/\dr\cJ) @= \vf_{T*}(\dr\cF_T/\dr\cJ)
  \endCD$$
 By the Nine Lemma, the induced map $v$ is an isomorphism.  Hence, since
$h$ is a presentation of colength $\delta'$, so is $h_T|\cJ$.

Since $\cCok(h_T)$ has support in $Y_T$, the map $h_T|\cJ$ is an
isomorphism off $Y_T$.  Also $\vf_T$ restricts to an isomorphism from
$\dr C_T-\dr Y_T$ onto $C_T-Y_T$.  Hence the points of $\dr C_T-\dr Y_T$
where ${\dr\cJ}^\star$ is not invertible are just those lying over the
points of $C_T-Y_T$ where $\cJ^\star$ is not invertible.  On the other
hand, since $\cI^\star$ is invertible along $Y_T$, its pullback
$\vf_T^*(\cI^\star)$ is invertible along $\dr Y_T$.  It follows that
$\cI^\star\ox(h_T|\cJ)^\star$ is a copresentation of colength $\delta'$.
Thus $\lambda_{\vf,h}$ is well defined.

 Note in passing that we have established the equation,
        $$V^m_{\vf,\dr\cF}=(1\x \theta_h)^{-1}W^{m,1}_{Y,\cF}.\eqno\Cs3.1)$$
 In particular,  $1\x\theta_h$ carries  $V^m_{\vf,\dr\cF}$ into
$W^{m,1}_{Y,\cF}$.

 \thm4 Let $C/S$ be a flat projective family of integral curves of
arithmetic genus $g$.  Let $\vf\:\dr C\to C$ be a birational $S$-map,
$\delta$ its genus change, $Y\subset C$ and $\dr Y\subset \dr C$ its
conductor subschemes.  Let $\cF$ and $\dr\cF$ be torsion-free rank-$1$
sheaves on $C/S$ and $\dr C/S$ of degrees $f$ and $f'$.  Assume
there is a presentation
        $h\:\cF\to\vf_{*}\dr\cF$
 whose adjoint $h^\sharp\:\vf^*\cF\to\dr\cF$ is surjective, and let
$h^\star\:\vf_*{\dr\cF}^\star\to\cF^\star$ be the associated map.
Fix $m$.  Set
  $$r:=m-1+f-2g+2,\ r':=m-1+f'-2g+2,\hbox{ and } \delta':=\delta+f'-f.$$

Then the following two diagrams are commutative:
 $$\CD
  V^m_{\vf,\dr\cF}               @>\lambda_{\vf,h}>> P_\vf^{r',\delta'} \\
    @V\rho\x1VV                                          @VV\pi_\vf V   \\
  W^{m-2\delta,1}_{\dr Y,\dr\cF} @>\alpha_{\dr\cF}>> \CJ^{r'}_{\dr C}
  \endCD\qquad\qquad\CD
  V^m_{\vf,\dr\cF}           @>\lambda_{\vf,h}>> P_\vf^{r',\delta'} \\
  @V1\x\theta_hVV                                      @VV\kappa_\vf V  \\
  W^{m,1}_{Y,\cF}             @>\alpha_\cF>>               \CJ^r_C
  \endCD$$
 The second is a product square if also $\cIm(h^\star)=\cI_Y\cF^\star$.
 \pf
 For each diagram, commutativity results directly from the definitions.
So there is a natural map,
        $$\gamma\:\  V^m_{\vf,\dr\cF}\to W^{m,1}_{Y,\cF}\x_{\CJ^r_C}
        P_\vf^{r',\delta'}.$$
 Assume that $\cIm(h^\star)= \cI_Y\cF^\star$.  Then we have to check
that $\gamma$ is an isomorphism.  To do so,
 we may change the base to an \'etale covering of $S$.  Thus we may
assume that the smooth locus of $C/S$ admits a section.  (Alternatively,
we could work directly, over the given $S$, with the map of presheaves
of $T$-points that we used to define $\gamma$.)

  First, let's check that $\gamma$ is surjective on $T$-points.  So take
 $S$-maps,
        $$(t_1, t_2)\:T\to\CJ^m_{C,Y}\x\Quot^1_\cF\and
        t\:T\to P_\vf^{r',\delta'},$$
 where $T$ is an arbitrary locally Noetherian $S$-scheme, such that
$(t_1, t_2)$ factors through $W^{m,1}_{Y,\cF}$ and
        $$\alpha_\cF(t_1, t_2)=\kappa_\vf t.\eqno\Cs4.1)$$
 Say  $\cJ\subset\cF_T$ defines $t_2$.  Since the smooth locus of
$C/S$ admits a section, there is a torsion-free rank-1 sheaf $\cI$ on
$C_T/T$ defining $t_1$, and there is a presentation $u\:\cK\to
\vf_{T*}\dr\cK$ defining $t$.  Since $(t_1, t_2)$ factors through
$W^{m,1}_{Y,\cF}$, the sheaf $\cI^\star$ is invertible along $Y_T$ and
where $\cJ^\star$ is not invertible.  Let $U\subset C_T$ be the open
subset where $\cI^\star$ is invertible.

 Because of \Cs4.1),
 the sheaves $\cK$ and $(\cI^\star\ox\cJ^\star)^\star$ differ by tensor
product with the pullback of an invertible sheaf on $T$.  We may replace
$\cI$ by its tensor product with this pullback, and so assume that
$\cK^\star=\cI^\star\ox\cJ^\star$.  Then $u$ induces a copresentation
        $$u^\star\:\vf_{T*}{\dr\cK}^\star\to\cI^\star\ox\cJ^\star.$$

 Define a subsheaf $\cL\subset\cJ^\star$ as follows.  On $C_T-Y_T$ define
$\cL$ to be $\cJ^\star$.  On $U$ define $\cL$ to be
$\cIm(u^\star)|U\ox(\cI^\star|U)^{-1}$.  The sheaf $\cL$ is defined
everywhere because $U\supset Y_T$.  Since $\cIm(u^\star)$ is equal to
$\cI^\star\ox\cJ^\star$ away from $Y_T$, the sheaf $\cL$ is well defined
and
        $$\cI^\star\ox\cL=\cIm(u^\star).\eqno\Cs4.2)$$
 Since $\cI^\star$ is invertible along $Y_T$, it follows
that there is a torsion-free rank-1 sheaf $\dr\cJ$
on $\dr C/T$ such that $\cL=\vf_{T*}{\dr\cJ}^\star$.

 Since $u^\star$ is a copresentation with colength $\delta'$, it follows
that the quotient $\cJ^\star/\cL$ is an $\cO_{Y_T}$-module of relative
length $\delta'$.  In particular, $\cL\supset\cI_{Y_T}\cJ^\star$.  Now,
the inclusion $\cJ\subset\cF_T$ induces a map $\cI_{Y_T}\cF_T^\star\to
\cI_{Y_T}\cJ^\star$.  Finally, our hypothesis that $\cIm(h^\star)=
\cI_Y\cF^\star$ implies that
 $ h_T^\star(\vf_{T*}{\dr\cF_T}^\star)=
\cI_{Y_T}\cF_T^\star$.
  Put together, these relations yield a map
$\vf_{T*}{\dr\cF_T}^\star\to\cL$.

Since $\cL=\vf_{T*}{\dr\cJ}^\star$, the map
$\vf_{T*}{\dr\cF_T}^\star\to\cL$ and the inclusion $\cL\into\cJ^\star$
yield maps $w$ and $v$, which help form the following commutative
diagram on $C_T$:
        $$\CD
        \cJ @>v>> \vf_{T*}\dr\cJ\\
        @VVV      @VVwV\\
        \cF_T @>h_T>>\vf_{T*}\dr\cF_T
        \endCD\eqno\Cs4.3)$$

By Lemma~(3.4), there is a unique map $\dr w\:\dr\cJ\to\dr\cF_T$ such
that $w=\vf_{T*}\dr w$.
 Since the relative colengths of $v$ and $h_T$ are both equal to
$\delta'$, that of $\dr w$ is 1.
  Hence $\cIm(\dr w)$ defines an $S$-map $\dr t_2\: T\to
\Quot^1_{\dr\cF}$.  It is clear from \Cs4.3) that $t_2=\theta_h\dr t_2$.
 Now,  $(t_1, t_2)$ factors through $W^{m,1}_{Y,\cF}$.  Hence, by \Cs3.1),
        $$(t_1,\dr t_2)\:T\to\CJ^m_{C,Y}\x\Quot^1_{\dr\cF}$$
 factors through $V^m_{\vf,\dr\cF}$.  Finally, it is clear from \Cs4.2)
that $t=\lambda_{\vf,h}(t_1,\dr t_2)$.  Hence,
        $$\gamma(t_1,\dr t_2)=((t_1, t_2),t).$$
 Thus $\gamma$ is indeed surjective on $T$-points.

To show that $\gamma$ is injective on $T$-points, let $t'_2\:T\to
\Quot^1_{\dr\cF}$ be another $S$-map such that $(t_1, t'_2)$ factors
through $V^m_{\vf,\dr\cF}$ and such that
        $$\gamma(t_1,t'_2)=((t_1, t_2),t).\eqno\Cs4.4)$$
 We have to show that $t'_2=\dr t_2$.  To do so, say $t'_2$ corresponds
to $\dr\cN\subset\dr\cF_T$, and set
        $$\cN:=h_T^{-1}(\vf_{T*}\dr\cN)\subset\cF_T.$$
  By definition, $\cN$ corresponds to $\theta_h(t'_2)$.  However,
$\theta_h(t'_2)=t_2$ by \Cs4.4).  Hence $\cN=\cJ$.

By definition, $\cI^\star\ox(h_T|\cN)^\star$ is the copresentation
corresponding to $\lambda_{\vf,h}(t_1,t'_2)$.  The latter is equal to
$t$ by \Cs4.4).  However, $t=\lambda_{\vf,h}(t_1,\dr t_2)$.  Hence,
$\cI^\star\ox(h_T|\cN)^\star$ and $u^\star$ are equivalent
copresentations.  Now, $\cL=\vf_{T*}{\dr\cJ}^\star$; so \Cs4.2) yields
$\cIm(u^\star)=\cI^\star\ox\vf_{T*}{\dr\cJ}^\star$.  Also $\cN=\cJ$.
Therefore, $\vf_{T*}{\dr\cJ}^\star$ and $\vf_{T*}{\dr\cN}^\star$ have
the same image in
 $\cJ^\star$.
  Hence, $w(\vf_{T*}{\dr\cJ})$ and
$\vf_{T*}\dr\cN$ are equal in $\vf_{T*}\dr\cF_T$.  Hence, $\dr w(\dr\cJ)$
and $\dr\cN$ are equal in $\dr\cF_T$.  So $t'_2=\dr t_2$.  Thus $\gamma$
is injective on $T$-points, and the proof is complete.

 \cor5 Let $C/S$ be a flat projective family of integral curves of
arithmetic genus $g$.  Let $\vf\:\dr C\to C$ be a birational $S$-map,
$\delta$ its genus change, $Y\subset C$ and $\dr Y\subset \dr C$ its
conductor subschemes, and $\vf^\sharp\:\cO_C\to\vf_*\cO_{\dr C}$ its
comorphism.  Fix $m$.  Assume that $C/S$ and $\dr C/S$ are Gorenstein,
set $V^m_{\vf,\dr C}:=V^m_{\vf,\cO_{\dr C}}$, and define the map
        $$\Lambda_\vf\:V^m_{\vf,\dr C}\to P_\vf^{m-1,\delta}
        \hbox{ by }(\cI,\dr\cJ)\mapsto\cI\ox(\vf^\sharp_T|\cJ)
 \hbox{ where } \cJ:=(\vf^\sharp_T)^{-1}(\vf_{T*}\dr\cJ).$$

Then the following two diagrams are commutative, and the second is a
product square:
        $$\CD
        V^m_{\vf,\dr C} @>\Lambda_\vf>> P_\vf^{m-1,\delta}\\
        @V\vf^*\x1VV                    @V\pi_\vf VV\\
        W^{m,1}_{\dr Y,\dr C} @>A_{\dr C}>> \CJ^{m-1}_{\dr C}
        \endCD\qquad\qquad\CD
     V^m_{\vf,\dr C} @>\Lambda_\vf>> P_\vf^{m-1,\delta}\\
         @V1\x\vf VV   @(\fiberbox)           @V\kappa_\vf VV\\
        W^{m,1}_{Y,C} @>A_C>> \CJ^{m-1}_C
        \endCD $$
 where $A_{\dr C}$ and $A_C$ are given by simple tensor product.
  \pf Observe that $\theta_{\vf^\sharp}=\vf$ under the natural
identifications,
        $$\Quot^1_{\cO_C}=C\and\Quot^1_{\cO_{\dr C}}=\dr C;$$
 so in the second diagram, $1\x\vf$ is well defined.  Since $C/S$ is
Gorenstein, a torsion-free rank-1 sheaf $\cI$ on $C$ is invertible
precisely where $\cI^\star$ is; similarly, for $\dr C/S$.  It follows
that $\vf^*\x1$, $A_{\dr C}$, and $A_C$ are well defined.  By the same
token, $\Lambda_\vf$ is well defined, because $\vf^\sharp_T|\cJ$ is a
presentation of colength $\delta$, as we noted toward the end of \Cs3).
Thus the diagrams make sense.  For each one, commutativity results
directly from the definitions.

To prove that the second diagram is a product square, apply
Theorem~\Cs4) with $\cF:=\cO_C$.  We have $\cIm({\vf^\sharp}^\star)=
\cI_Y\bigomega_C$ by (3.1.1) since $C/S$ is Gorenstein.  Now, since the
dualizing sheaf
 $\bigomega_C$
 is invertible, tensoring with it defines an
isomorphism,
        $$\tau'\:P_\vf^{m-1-2g+2,\delta}\risom P_\vf^{m-1,\delta}.$$
 We have $\Lambda_\vf=\tau'\circ\lambda_{\vf,\vf^\sharp}$; indeed, we
can prove it the same way that we proved that
 $A_C=\tau\circ\alpha_C$
 in the proof of
Corollary~(2.6).  Hence the second diagram is a product square because,
in Theorem~\Cs4), the second diagram is one.  The proof is now complete.

\sct6 Double points

\prp1 Let $C$ be a projective integral curve over an algebraically
closed field $k$, and $Q$ a singular point.  Then the following
conditions are equivalent:\smallbreak
 \item i the point $Q$ is a double point;
 \item ii there exist variables $u,v$ and formal power series
$a(u),\ b(u)$ such that
        $$\wh\cO_{C,Q}=k[[u,v]]/(v^2+a(u)uv+b(u)u^2);$$
 \item iii given any torsion-free rank-$1$ sheaf $\cI$ on $C$, at most two
elements are required to generate its stalk  $\cI_Q;$
 \item iv we have $\dim_k\cO_{C,Q}/\Im_{C,Q}^3\le5;$
 \item v
 along the blowup $\phi\:\wt C\to C$ at $Q$, the genus drop
$\delta_\phi$ is equal to $1$.
 \pf For convenience, set $A:=\cO_{C,Q}$, set $\Im:=\Im_{C,Q}$,
set $I:=\cI_Q$, and set $B:=(\phi_*\cO_{\wt C})_Q$.

 Assume (i).  Since $k$ is infinite, there is by \cite{ZS60, Thm.~22,
\p.294} a parameter $u$ of $A$ such that the ideal $(u)$ has
multiplicity 2.  Since $u$ is not a zero-divisor, $A/(u)$ has length 2.
Moreover, the completion $\wh A$ contains the formal power series ring
$k[[u]]$.  Hence $\wh A$ is a free module of rank 2 over $k[[u]]$.
Therefore, $\wh A$ is of the form $k[[u,v]]/(F)$ where $F$ is a monic
polynomial of degree 2 in $v$ with coefficients in $k[[u]]$.  Since $A$
is singular, $F$ must lie in $(u,v)^2$; so $F$ must be of the form
indicated in (ii).  Thus (ii) holds.

Assume (ii) and consider (iii).  Then the completion $\wh I$ is a free
$k[[u]]$-module.  It is of rank 2 because $\wh A$ is.  Thus (iii) holds.

Assume (iii).  Then (iv) holds because we always have
     $$\dim_kA/\Im^3=\dim_kA/\Im+
        \dim_k\Im/\Im^2+\dim_k\Im^2/\Im^3.\eqno\Cs1.1)$$

Assume (iv).  Then the right side of \Cs1.1) is at most 5.  Now, since
$A$ is singular, then both $\Im/\Im^2$ and $\Im^2/\Im^3$ have dimension
at least 2, so exactly 2.  Hence, $\wh A$ is of the form
$k[[u,v]]/(F)$ where $F\in(u,v)^2$, but $F\notin(u,v)^3$.  Consequently,
$\wh A$ has multiplicity 2.  So (i) holds.

Assume (ii) and consider (v).  Set $w:=v/u$.
  Then $\wh
B=\wh A[w]$.  Now, $\wh A[w]=\wh A+kw$ since $w^2+a(u)w+b(u)=0$.
Hence, (v) holds.

Assume (v).  Then we have
        $$\dim_kB/\Im=\dim_kB/A+\dim_kA/\Im=1+1=2.$$
 Now, $\dim_kB/\Im B>1$; otherwise, $A/\Im\to B/\Im B$ is surjective,
and so, by Nakayama's lemma, $A\to B$ is surjective, contrary to
our hypothesis.  Hence $\Im=\Im B=tB$ for some $t$.  Therefore,
        $$B/\Im\risom \Im/\Im^2\risom \Im^2/\Im^3\risom\cdots.$$
 Consequently, Equation \Cs1.1) yields (iv).  Alternatively, (i) holds
by the definition of multiplicity.  The proof of the proposition is now
complete.

 \prp2 Let $C/S$ be a flat projective family of integral curves, and $Q$
a flat closed subscheme of the singular locus of $C/S$ {\rm(}which is
the subscheme of $C$ defined by the first Fitting ideal of the sheaf of
differentials $\Omega^1_{C/S}${\rm)}.  Assume that $Q(s)$ is a double
point of $C(s)$ for each geometric point $s$ of $S$.  Let
$\cM\subset\cO_C$ denote the ideal of $Q$; set $\cB:=\cHom(\cM,\cM)$ and
$C^*:=\Spec(\cB)$.  Let $\wt C$ be the blowup of $C$ along $Q$.  Then
$C^*$ and $\wt C$ are $C$-isomorphic, flat projective families of
integral curves with genus drop $1$ from $C$; moreover, forming $C^*$,
$\wt C$, and the isomorphism commutes with changing $S$.
 \pf
  To prove that $\wt C/S$ is flat and that, given a base change $T\to
S$, the induced map $(C\x T)\wt{\phantom{c}}\to\wt C\x T$ is an
isomorphism, we may replace $S$ by $\Spec\cO_{S,s}$ where $s$ is an
arbitrary point of $S$.  Then, we may replace $\cO_{S,s}$ by a flat
local over-ring whose residue field is algebraically closed; such a ring
exists by \cite{EGA, O$_{\rm III}$ 10.3.1, \p.20}.  We may now embed $C$
into a projective family of surfaces $F/S$ that is smooth along $Q$
since $Q(s)$ is a surficial singularity.  Indeed, we can construct $F/S$
using Bertini's theorem on the hypersurfaces of high degree containing
$C$ in an ambient projective space; for more details, see the beginning
of the proof of \cite{EGK99, (3.4)}.

Consider the nested sequence $Q\subset C\subset F$, and let $\cI_{Q/F}$
and $\cI_{C/F}$ denote the ideals.  Restrict
to $Q$ the usual right exact sequence, getting
 $$\cI_{C/F}/\cI_{C/F}^2|Q\to\Omega^1_{F/S}|Q\to\Omega^1_{C/S}|Q\to0.$$
 The middle term is locally free of rank 2; in fact, it is equal to
$\cI_{Q/F}/\cI_{Q/F}^2$.  So we can use this sequence to compute the
first Fitting ideal of the third term.  We find that, since $Q$ lies in
the singular locus of $C/S$, we have $\cI_{C/F}\subset\cI_{Q/F}^2$.
Moreover, $\cI_{C/F}(s)\not\subset\cI_{Q/F}^3(s)$ for each $s\in S$
because $Q(s)$ is a double point.

Form the blowup map $\beta\:\wt{F}\to F$ along $Q$, and denote the
exceptional divisor by $E$.  Now, $\beta^{-1}C$ is a relative effective
Cartier divisor on a neighborhood of $E$ in $\wt{F}$ because
$F/S$ is smooth along $Q$.
Moreover, $\beta^{-1}C> 2E$ because $\cI_{C/F}\subset\cI_{Q/F}^2$.
Let $R$ be the residual scheme to $2E$ in $\beta^{-1}C$; so, $R =
\beta^{-1}C-2E$ on a neighborhood of $E$, and $R=\wt C$ off $E$.
Then $R$ is flat, and forming it commutes with changing $S$.
Now, $E$ restricts to a Cartier divisor on $R$ because
$\cI_{C/F}(s)\not\subset\cI_{Q/F}^3(s)$ for each $s\in S$. Since
$E$ restricts to a Cartier divisor on $\wt{C}$ as well, $\wt{C}=R$
because each side is equal to its closure in $\wt{F}$ of its restriction to
$\wt{F}-E$.
Hence $\wt{C}$ is flat, and forming it commutes with changing $S$.
Moreover, by Lemma~(6.1), the genus drop from $C$ to $\wt{C}$ is 1.

It remains to prove that $C^*$ is equal to $\wt{C}$.  Let $\phi\:\wt
C\to C$ denote the structure map.  Set $\wt{Q}:=\phi^{-1}Q$, and let
$\wt{\cM}$ denote its ideal.
 The genus drop of $\phi$ is 1, so $\wt{Q}/S$ is flat of relative length
2 by \cite{AK90, Prop.~2, \p.17}.  So $\cM=\phi_*\wt{\cM}$.  Hence
  $$\cB=\cHom(\phi_*\wt{\cM},\,\phi_*\wt{\cM}).$$
 By Lemma (3.4), or rather its local version, we have
  $$\cHom(\phi_*\wt{\cM},\,\phi_*\wt{\cM})=\phi_*Hom(\wt{\cM},\,\wt{\cM}).$$
 Now, $\wt{\cM}$ is invertible; so $Hom(\wt{\cM},\,\wt{\cM})=\cO_{\wt{C}}$.
Hence $\cB=\phi_*\cO_{\wt{C}}$.  Now, $\Spec(\phi_*\cO_{\wt{C}})=\wt{C}$
since $\phi$ is finite, so affine. Therefore, $C^*=\wt{C}$, and the
proof is complete.

 \thm3 Let $C/S$ be a flat projective family of integral curves, and $Q$
a flat closed subscheme of the singular locus of $C/S$.  Assume that
$Q(s)$ is a double point of $C(s)$ for each geometric point $s$ of $S$.
  Let $\phi\:\wt C\to C$ be the blowup along $Q$.  If the smooth locus
of $C/S$ admits a section, then the presentation scheme $P^{m,1}_\phi$
is a $\IP^1$-bundle over the compactified Jacobian $\CJ^m_{\wt C}$.
 \pf By  Proposition~\Cs2),
$\phi$ is a birational $S$-map with genus drop 1.  So
we may consider its presentation scheme $P^{m,1}_\phi$.
Assume that the smooth locus of $C/S$ admits a section.
Then $\wt C\x\CJ^m_{\wt C}$ admits
a universal sheaf $\wt\cI$.  Set
        $$\cF:=((\phi\x1)_*\wt\cI)\big|(Q\x\CJ^m_{\wt C}).$$
 Then $P^{m,1}_\phi=\IP(\cF)$ by \cite{AK90, Thm.~10, \p.21}, because
the genus drop
 along
 $\phi$ is 1.  We have to prove
that $\cF$ is locally free of rank 2.

Set $\wt Q:=\phi^{-1}Q$, and denote its ideal by $\wt\cM$.  Since
$\phi$ is a birational $S$-map with genus drop 1, by
\cite{AK90, Prop.~2, \p.17} the scheme $\wt Q$ is flat.
In addition, $\wt\cM$ is invertible since $\phi$ is the blowup map at $Q$.
Consider the sequence,
        $$0\to p^*\wt\cM\ox\wt\cI\to\wt\cI\to
        \wt\cI\bigm|(\wt Q\x\CJ^m_{\wt C})\to 0,$$
 where $p\:\wt C\x\CJ^m_{\wt C}\to\wt C$ is the projection.  Since
$\wt\cM$ is invertible,
 and forming it commutes with changing $S$, and since
 $\wt\cI$ is torsion-free rank-1 on $\wt C\x\CJ^m_{\wt
C}/\CJ^m_{\wt C}$, the sequence is exact on the fibers over $\CJ^m_{\wt
C}$. Hence (it is exact, and) the third term is flat
over $\CJ^m_{\wt C}$ by \cite{EGA, IV$_3$ 11.3.7, \p.135}.  So $\cF$ is
locally free.

To compute the rank of $\cF$, we may assume that $S$ is the spectrum
of an algebraically closed field.
The rank of $\cF$ is the minimal number of generators at $Q$ of any
fiber of $(\phi\x1)_*\wt\cI$ over $\CJ^m_{\wt C}$.  By Lemma~(3.8),
no fiber is invertible at $Q$, hence the rank is at least 2.
On the other hand, the rank is at most 2 by
Condition (iii) of Proposition~\Cs1).  Thus the rank is 2,
and the proof is complete.

 \lem4 Let $C$ be a projective integral curve over an algebraically
closed field, $Q$ a double point of $C$, and $\cI$ a torsion-free rank-$1$
sheaf on $C$.  Denote by $\phi\:\wt C\to C$ the blowup at $Q$.  Then the
following assertions hold.
 \part 1 Let $\vf\:\dr C\to C$ be a birational map that is nontrivial over
$Q$.  Then there are three possibilities for the fiber $\vf^{-1}Q$: two
simple points, one simple point, or one double point.  Moreover, $\vf$
factors through the blowup $\phi\:\wt C\to C$.
 \part 2 There exists a torsion-free rank-$1$ sheaf $\wt\cI$ on $\wt C$
such that $\phi_*\wt\cI=\cI$ if and only if $\cI$ is not invertible at $Q$.
If $\wt\cI$ exists, then it is unique.
 \part 3 There exist a unique birational map $\vf\:\dr C\to C$ that is
trivial off $Q$, and a unique torsion-free rank-$1$ sheaf $\dr\cI$ on $\dr
C$ such that $\dr\cI$ is invertible along $\vf^{-1}Q$ and $\vf_*\dr\cI=\cI$.
 \pf
 To prove (1), form the sum of the multiplicities of the points in
$\vf^{-1}Q$.  This sum is always at most the multiplicity of $Q$ (look
at the pullback of an element $u\in \cO_{C,Q}$ such that the ideal
$u\cO_{C,Q}$ is of the same multiplicity as $Q$).  The first assertion
follows immediately.

Consider the second assertion of (1) and also (2).  Proposition~\Cs2)
says that $\wt C$ is isomorphic over $C$ to $C^*$, where
$C^*:=\Spec(\cHom(\cM,\cM))$.  So the second assertion of (1) follows
from Part~(3) of Lemma~(3.7), and (2) follows from Lemma~(3.8).

Consider (3).  Suppose $\cI$ is invertible at $Q$.  Then $1_C$ is a
suitable $\vf$.  There can be no other; indeed, any other would factor
through $\wt C$ by (1), and so lead to a contradiction of (2).
So suppose $\cI$ is not invertible at $Q$.  Then by (2) there is an
$\wt\cI$ on $\wt C$ such that $\phi_*\wt\cI=\cI$.

Let $\delta$ be the genus drop to the (partial) normalization of $C$ at
$Q$, and proceed by induction on $\delta$.  Suppose $\delta=1$.  Then
$\phi\:\wt C\to C$ is the only nontrivial birational map trivial off
$Q$, and $\wt\cI$ is invertible along $\phi^{-1}Q$.  So (3) holds in
this case.

Suppose $\delta\ge2$.  By Condition~(v) of Proposition~\Cs1), the genus
drop along $\phi$ is 1.  So (1) implies that $\phi^{-1}Q$ consists of
one double point.  Hence, by induction, we may assume that (3) holds for
$\wt C$ and $\wt\cI$.  The existence of a $\vf\:\dr C\to C$ and of an
$\dr\cI$ follows directly, and their uniqueness follows because any
suitable $\vf$ must factor through an $\vf'\:\dr C\to\wt C$ by (1), and
any suitable $\dr\cI$ must satisfy $\vf'_*\dr\cI=\wt\cI$ by Lemma~(3.4).
Thus (3) holds, and the proof is complete.

 \lem5 Let $C$ be a projective integral curve over an algebraically
closed field, and $Y\subset C$ a closed subset.  Fix $n$.  Then the
following assertions hold.
 \part1 Let $\vf\:\dr C\to C$ be a birational map that is
trivial off $Y$.  Then
        $$\epsilon_\vf(\CJ^{n-\delta_\vf}_{\dr C,\vf^{-1}Y})
        \subset\CJ^{n,\delta_\vf}_{C,Y}-\CJ^{n,\delta_\vf-1}_{C,Y}.$$
 \part2
 Let $\vf_1\:\dr C_1\to C$ and $\vf_2\:\dr C_2\to C$ be birational maps
that are trivial off $Y$.  Then there exists an isomorphism $\mu\:\dr
C_1\risom\dr C_2$ such that $\vf_1=\vf_2\circ\mu$ if and only if
        $$\epsilon_{\vf_1}(\CJ^{n-\delta_{\vf_1}}_{\dr C_1,\vf_1^{-1}Y})
        \cap\epsilon_{\vf_2}(\CJ^{n-\delta_{\vf_2}}_{\dr C_2,\vf_2^{-1}Y})
        \neq\emptyset.\eqno\Cs5.1)$$
 \pf Let's prove (1).
  Let $t\in \CJ^{n-\delta_\vf}_{\dr C,\vf^{-1}Y}$ be a closed point,
and $\dr\cK$ a corresponding sheaf on $\dr C$.  Set $\cK:=\vf_*\dr\cK$.
Since $\vf$ is trivial off $Y$, its conductor subscheme on $\dr C$ is
 set-theoretically
 contained in $\vf^{-1}Y$.  So ${\dr\cK}^\star$ is
invertible along the conductor subscheme.  Then Part~(1) of
Proposition~(4.3) says that $\epsilon_\vf(t)\in\CJ^{n,k}_{C,Y}$ if and
only if $k\ge\delta_\vf$.
  Thus (1) is proved.

 Let's prove (2).  If an isomorphism $\mu$ exists, then certainly
\Cs5.1) holds; in fact, the two images coincide.

Conversely, suppose \Cs5.1) holds.  Then there exist three
torsion-free rank-1 sheaves, $\cI$ on $C$ and
 $\dr\cI_1$ on $\dr C_1$ and
$\dr\cI_2$ on $\dr C_2$,
 such that $\dr\cI_i$ is invertible along $\vf_i^{-1}Y$
and $\cI=\vf_{i*}\dr\cI_i$.  Write $\cI=\cJ\ox\cL$ where $\cJ$ is the
ideal of a closed subscheme $Z\subset C$ with support in $Y$ and $\cL$
is a torsion-free rank-1 sheaf on $C$ that is invertible along $Y$. Set
$\dr\cJ_i:=\cJ\cO_{\dr C_i}$.  Since $\cL$ is invertible along $Y$, the
pullback $\vf_i^*\cL$ is a torsion-free rank-1 sheaf on $\dr C_i$ that
is invertible along $\vf_i^{-1}Y$.  Hence $\dr\cI_i=
\dr\cJ_i\ox\vf_i^*\cL$.  Since $\dr\cI_i$ is invertible along
$\vf_i^{-1}Y$, so is $\dr\cJ_i$.  Moreover, since
$\cI=\vf_{i*}\dr\cI_i$, the projection formula yields
$\cJ=\vf_{i*}\dr\cJ_i$.

Let $\phi\:\wt C\to C$ be the blowup of $C$ along $Z$.  By the universal
property, since $\dr\cJ_i$ is invertible, $\vf_i$ factors through
$\phi$.  Let $\psi_i\:\dr C_i\to\wt C$ denote the induced map.  Set
$\wt\cJ:=\cJ\cO_{\wt C}$.  Then $\wt\cJ=\psi_{i*}\dr\cJ_i$ since
$\cJ=\vf_{i*}\dr\cJ_i$.  Since $\wt\cJ$ and $\dr\cJ_i$ are invertible,
therefore the comorphism
 $\psi_i^\sharp:\cO_{\wt
C}\to\psi_{i*}\cO_{\dr C_i}$
 is an isomorphism.  So $\psi_i$ is an
isomorphism.  Let $\mu:=\psi_2^{-1}\circ\psi_1$.  Then $\mu\:\dr
C_1\to\dr C_2$ is the desired isomorphism.  The proof is now complete.

 \prp6 Let $C$ be a projective integral curve over an algebraically
 closed field, $Y\subset C$ a closed subset.  Assume that $Y$ consists
of points of multiplicity at most $2$.  Fix $k\geq 0$.  Then
        $$\CJ^{n,k}_{C,Y}=\coprod_{\{\vf|\delta_\vf\le k\}}
        \epsilon_\vf(\CJ^{n-\delta_\vf}_{\dr C,\vf^{-1}Y})$$
 where the disjoint union runs over all birational maps
 $\vf\:\dr C\to C$ that have genus drop $\delta_\vf$ at most $k$ and that
are trivial off $Y$.
 \pf
 By convention, $\CJ^{n,-1}_{C,Y}=\emptyset$; see (4.1).  By
Part~(1) of Proposition~(4.4),
        $$\CJ^n_C=\coprod_{k\ge0}(\CJ^{n,k}_{C,Y}-\CJ^{n,k-1}_{C,Y}).$$
 By Parts~(1)~and~(2) of Lemma~\Cs5), we have only to prove that
        $$\CJ^n_C=\bigcup_{\{\vf|\delta_\vf\ge 0\}}
        \epsilon_\vf(\CJ^{n-\delta_\vf}_{\dr C,\vf^{-1}Y}),\eqno\Cs6.1)$$
 where the union runs over all birational maps
$\vf\:\dr C\to C$ that are trivial off $Y$.

Let $\delta_Y$ denote the genus drop
 to
 the (partial) normalization of $C$ along $Y$.  We'll prove \Cs6.1) by
induction on $\delta_Y$. If
 $\delta_Y=0$,
 then $C$ is smooth along $Y$.
So $1_C$ is the only birational map $\dr C\to C$ that is trivial off
$Y$.  Moreover, every
 torsion-free
 rank-1 sheaf on $C$ is invertible
along $Y$.  Hence \Cs6.1) holds in this case.

Suppose $\delta_Y>0$.  Let $t$ be a
 closed point
 of $\CJ^n_C$,
and $\cK$
 a
 corresponding sheaf.  We have to construct a
birational map $\vf\:\dr C\to C$ such that
 $$t\in\epsilon_\vf(\CJ^{n-\delta_\vf}_{\dr C,\vf^{-1}Y}).\eqno\Cs6.2)$$

 Suppose that $\cK$ is invertible at every point of $Y$.  Take
$\vf:=1_C$.  Then $\delta_\vf=0$.  Now, $C$ is Gorenstein along $Y$ by
Proposition~\Cs1).  Hence $\cK^\star$ too is invertible along $Y$.
Hence \Cs6.2) holds.

Suppose that $\cK$ is not invertible at $Q\in Y$.
 Then $Q$ is double.  By Part~(3) of Lemma~\Cs4), there exist a
birational map $\phi\:C'\to C$ that is trivial off $Q$, and a
 torsion-free
 rank-1 sheaf $\cK'$ on $C'$ such that $\cK'$ is invertible
along $\phi^{-1}Q$ and $\phi_*\cK'=\cK$.
 Set $Y':=\phi^{-1}Y$.
 By Part~(1) of Lemma~\Cs4), the points of $Y'$ have multiplicity at
most 2 in $C'$.  Since $\cK$ is not invertible at $Q$, the genus drop
along $\phi$ is positive.  Hence, the genus drop to the normalization of
$C'$ along $Y'$ is strictly less than $\delta_Y$.
 \par

 By induction, \Cs6.1) holds for $C'$ and $Y'$.  So there exist a
birational map $\psi\:\dr C\to C'$ that is trivial off $Y'$, and a
torsion-free rank-1 sheaf $\dr\cK$ on $\dr C$ such that ${\dr\cK}^\star$
is invertible along $\psi^{-1}Y'$ and $\psi_*\dr\cK=\cK'$.  Set
$\vf:=\phi\circ\psi$.  Then ${\dr\cK}^\star$ is invertible along
$\vf^{-1}Y$ and $\vf_*\dr\cK=\cK$.  So \Cs6.2) holds.  Thus \Cs6.1)
holds, and the proof is complete.

 \thm7 Let $C$ be a projective integral curve of arithmetic genus $g$
over an algebraically closed field.  Fix $n$.  Assume that $C$ is
Gorenstein.  Then $C$ has double points at worst if and only if
$\dim(\CJ^n_C-\CJ^{n,1}_{C,C})\le g-2$.  If so, then equality holds,
unless $C$ is smooth or has just one singular point, which is either an
ordinary double point or cusp; in the exceptional case,
$\CJ^n_C=\CJ^{n,1}_{C,C}$.
 \pf
 First suppose $\dim(\CJ^n_C-\CJ^{n,1}_{C,C})\le g-2$.
 Let $Q\in C$ be a singular point.  By Parts~(2) and (3) of Lemma~(3.7),
 there exists a unique birational map $\vf\:\dr C\to C$
 that is trivial off $Q$ and has $\delta_\vf=1$.
 The associated map $\epsilon_\vf$
embeds the Jacobian $J_{\dr C}^{n-1}$ in $\CJ_C^n$ as a locally closed
subset of dimension $g-1$.  Hence this subset meets $\CJ^{n,1}_{C,C}$.
 So there exist an invertible sheaf $\dr\cL$ on $\dr C$, a torsion-free
rank-1 sheaf $\cL$ on $C$, and a subsheaf $\cI\subset\bigomega_C$ of
colength 1 such that $\cL^\star$ is invertible and
$\vf_*\dr\cL\cong(\cL^\star\ox\cI^\star)^\star$.

  As $\bigomega_C$ is invertible, so is $\cL$.
 Set $\cM:=\cI\ox\bigomega_C^{-1}$.  Then $\cM$ is the ideal of a point
of $C$, and Equations (2.6.3) and (2.6.4) yield the equation,
        $$(\cL^\star\ox\cI^\star)^\star=\cL\ox\cM.\eqno\Cs7.1)$$
 By the projection formula, $\cM\cong\vf_*(\dr\cL\ox\vf^*\cL^{-1})$.
So $\cM$ is contained in the conductor ideal of $\vf$, whence must be
equal to it because $\vf$ is nontrivial.   Hence
 $\cM$ is the ideal of $Q$.  Moreover, $\vf^*\cM$ has an invertible
quotient, namely, $\dr\cL\ox\vf^*\cL^{-1}$.  So $\vf$ factors through
the blowup $\phi$ of $C$ at $Q$; hence, $\vf=\phi$ because
$\delta_\vf=1$.  By Proposition~\Cs1), the point $Q$ is double.

 Conversely, suppose that every singular point $Q$ of $C$ is double.
Then Proposition~\Cs6) implies that
        $$\CJ_C^n-\CJ^{n,1}_{C,C}=\bigcup_{\{\vf|\delta_\vf\ge 2\}}
  \epsilon_\vf\circ\iota(J_{\dr C}^{2g-2-n+\delta_\vf}).\eqno\Cs7.2)$$
 Hence, either $\CJ_C^n=\CJ^{n,1}_{C,C}$ or else there exists a $\vf$
with $\delta_\vf\ge 2$.

Suppose that $\vf$ has $\delta_\vf\ge 2$.  Then $\vf$ factors through
the blowup $\phi\:\wt C\to C$ at $Q$ by Proposition~(6.2) and Part (3)
of Proposition~(3.7).  Moreover, $\delta_\phi=1$ by Proposition~(6.1),
and $\phi^{-1}Q$ is a double point by Part (1) of Lemma (6.4).  Form the
blowup $\phi'$ of $\wt C$ at $\phi^{-1}Q$.  Then $\delta_{\phi'}=1$ by
Proposition~(6.1).  Set
 $\vf':=\phi\circ\phi'$.
  Then $\delta_{\vf'}=2$.
 Now, $\epsilon_\vf\circ\iota$ is an embedding, and the dimension of
$J_{\dr C}^{2g-2-n+\delta_\vf}$ is $g-\delta_{\vf}$.  Therefore, \Cs7.2)
yields $\dim(\CJ^n_C-\CJ^{n,1}_{C,C})= g-2$.

On the other hand, suppose that $\CJ_C^n=\CJ^{n,1}_{C,C}$.  Then every
$\vf\:\dr C\to C$ has genus drop at most 1.  In particular, the
normalization map does.  Hence $C$ is smooth or has just one singular
point, either an ordinary double point or a cusp.  Thus the theorem is
proved.

 \cor8 Let $C$ be a projective integral curve of arithmetic genus $g$
over an algebraically closed field.  Fix $n$.  Assume that $C$ is
Gorenstein.  Let $V\subset\CJ^n_C$ be the image of the Abel
map of bidegree $(n+1,1)$,
        $$A_C\:J^{n+1}_C\x C \to \CJ^n_C, \hbox{ defined by }
        (\cL,\cM)\mapsto\cL\ox\cM.$$
 Then $C$ has double points at worst if and only if $\dim(\CJ^n_C-V)\le
g-2$.
 \pf Since the dualizing sheaf $\bigomega$ is invertible, sending a
point to the tensor product of its ideal with $\bigomega$, we get a
natural isomorphism,
        $$\sigma\:C\risom\Quot^1_\bigomega.$$
 Moreover, a torsion-free rank-1 sheaf $\cL$ is invertible if and only
if $\cL^\star$ is.  Hence,
 $1_{\CJ^{n+1}_C}\x\sigma$
 induces an
isomorphism,
        $$\chi\:J^{n+1}_C\x C\risom W^{n+1,1}_{C,\bigomega}.$$
 Equation~\Cs7.1) says in other words that $\alpha_\bigomega\circ\chi=
A_C$.  Hence $\CJ^{n,1}_{C,C}=V$.  Therefore, Theorem~\Cs7) yields the
assertion.

 \references
 \serial{adv}{Adv. Math.}
 \serial{ajm}{Amer. J. Math.}
 \serial{ansci}{Ann. Sci. \'Ecole Norm. Sup. (4)}
 \serial{bams}{Bull. Amer. Math. Soc.}
 \serial{bpm}{Birkh\"auser Prog. Math.}
 \serial{ca}{Comm. Alg.}
 \serial{compo}{Compositio Math.}
 \serial{India}{Proc. Indian Acad. Sci. Sect. A Math. Sci.}
 \def\GrothFest#1 #2 {\unskip, \rm
 in ``The Grothendieck Festschrift. Vol.~#1,'' P.~Cartier, L.~Illusie,
N.M.~Katz, G.~Laumon, Y.~Manin, K.A.~Ribet (eds.) \bpm 86 1990 #2 }
 \def\osloii{\unskip, \rm in ``Real and complex singularities. Proceedings,
Oslo 1976,'' P. Holm (ed.), Sijthoff \& Noordhoff, 1977, pp.~}
 
 \def\splnm #1 #2 {Springer Lecture Notes in Math. {\bf#1}\ (#2)}

AIK76
 {A. Altman, A. Iarrobino and S. Kleiman},
 Irreducibility of the compactified Jacobian
 \osloii 1--12.

AK76
 A. Altman and S. Kleiman,
 Compactifying the Jacobian,
 \bams 82(6) 1976 947--49

AK79II
 A. Altman and S. Kleiman,
 Compactifying the Picard scheme II,
 \ajm 101 1979 10--41

AK80
 A. Altman and S. Kleiman,
 Compactifying the Picard scheme,
 \adv 35 1980 50--112

AK90
 A. Altman and S. Kleiman,
 The presentation functor and the compactified Jacobian
 \GrothFest I 15--32

E95
 D. Eisenbud,
 Commutative Algebra,
 Springer GTM {\bf 150} (1995).

EGK99
 {E. Esteves, M. Gagn\'e and S. Kleiman},
 Autoduality of the compactified Jacobian,
 preprint.

EGA
 A. Grothendieck and J. Dieudonn\'e,
 El\'ements de G\'eom\'etrie Alg\'ebrique,
 Publ. Math. IHES, O$_{\rm III}$, Vol. 11, 1961;
III$_2$, Vol. 17, 1963; IV$_3$, Vol. 28, 1966; and
IV$_4$, Vol. 32, 1967.

KK81
 {S. Kleiman and H. Kleppe},
 Reducibility of the compactified Jacobian,
 \compo 43(2) 1981 277--280

R80
 C. Rego,
 The compactified Jacobian,
 \ansci 13(2) 1980 211-223

ZS60
 O. Zariski and P. Samuel,
 Commutative Algebra, Vol II,
 Van Nostrand, 1960.

\endreferences

\mark{E}\bye